\documentclass[12pt]{article}

\usepackage{verbatim}
\usepackage{amsmath}
\usepackage{amsthm}
\usepackage{amsfonts}
\usepackage{mathrsfs}
\usepackage{pxfonts}
\usepackage{amssymb}
\usepackage{amsbsy}
\usepackage{graphicx}
\usepackage{setspace}
\newtheorem{thm}{Theorem}[section]
\newtheorem{cor}[thm]{Corollary}
\newtheorem{lem}[thm]{Lemma}

\newtheorem{rem}[thm]{Remark}

\newtheorem{conj}[thm]{Conjecture}

\newcommand{\ou}{\omegaup}

\newcommand{\D}{\displaystyle}

\newcommand{\wc}{\overset{*}{\longrightarrow}}

\newcommand{\cj}[1]{\overline{#1}}

\newcommand{\ext}{\mathrm{ext}}

\newcommand{\ointcc}{\varointctrclockwise}

\newenvironment{demo}[1]{\renewcommand{\proofname}{#1}\begin{proof}}{\end{proof}}
\numberwithin{equation}{section}
\usepackage[left=1in,right=1in,top=1.3in,bottom=1in,footnotesep=18pt]{geometry}

\title{On the asymptotic behavior of Faber polynomials for domains with piecewise analytic boundary}
\author{Erwin Mi\~{n}a-D\'{\i}az}

\begin{document}
\maketitle

\abstract{For a function $\psi(w)$ analytic and
univalent in $\left\{w:1<|w|<\infty\right\}$ with
a simple pole at $\infty$ and a continuous
extension to $\{w:|w|\geq 1\}$, we consider the
Faber polynomials $F_n(z)$, $n=0,1,2,\ldots$,
associated to $\psi$ via their generating
function
$\psi'(w)/\left(\psi(w)-z\right)=\sum_{n=0}^\infty
F_n(z)w^{-(n+1)}$. Assuming that $\psi$ maps the
unit circle $\mathbb{T}_1$ onto a piecewise
analytic curve $L$ whose exterior domain has no
outward-pointing cusps, and under an additional
assumption concerning the ``Lehman expansion" of
$\psi$ about those points of $\mathbb{T}_1$
mapped onto corners of $L$, we obtain asymptotic
formulas for $F_n$ that yield fine results on the
location, limiting distribution and accumulation
points of the zeros of the Faber polynomials. The
asymptotic formulas are shown to hold uniformly
and the exact rate of decay of the error terms
involved is provided.\vspace{.3cm}}



\section{Introduction}
Let $\phi$ be a function with a Laurent expansion at
$\infty$ of the form
\begin{equation}\label{FabEq83}
\phi(z)=b_1z+b_0+\frac{b_{-1}}{z}+\frac{b_{-2}}{z^2}+\cdots,\quad
b_1\not=0, \quad
\limsup_{n\to\infty}|b_{-n}|^{1/n}<\infty.
\end{equation}
The $n$th Faber polynomial $F_n(z)$, $n=0,1,\ldots\,$,
associated with $\phi$ is the polynomial part of the
Laurent expansion at infinity of the function
$[\phi(z)]^n$.

We shall frequently use the following notation: given
$r\geq 0$,
\[
\mathbb{T}_r:=\{w:|w|=r\},\quad \Delta_r:=\{w:r<|w|\leq
\infty\}.
\]

The inverse function of $\phi$, denoted by $\psi$, is
well-defined in a neighborhood of $\infty$, and there is a
smallest number $\rho<\infty$ such that $\psi$ has an
analytic and \emph{univalent} continuation to
$\Delta_\rho\setminus\{\infty\}$. If $\rho=0$, then $\phi$
is linear and $F_n(z)=(b_1z+b_0)^n$. Being this case a
trivial one,
 we assume hereafter that $\phi$ has been normalized
so that $\rho=1$.

Then, the function $\psi$ maps $\Delta_1$ conformally onto
a simply-connected domain $\Omega$, and consequently,
$\phi$ has a conformal extension to $\Omega$, with
$\phi(\Omega)=\Delta_1$. Conversely, by the Riemann mapping
theorem, given any simply-connected neighborhood $\Omega$
of $\infty$ whose boundary contains more than one point,
there is, up to a multiplicative unimodular constant, a
unique conformal map $\phi$ of $\Omega$ onto $\Delta_1$
that complies with (\ref{FabEq83}). Hence, Faber
polynomials are often introduced as being generated by
simply-connected neighborhoods of $\infty$.

The function (in the variable $w$)
$\psi'(w)/\left(\psi(w)-z\right)$ is called the generating
function of the Faber polynomials, since as shown by Faber
\cite{Faber} (see also \cite{Ullman}), its Laurent
expansion at $\infty$ is
\begin{equation}\label{FabEq84}
\frac{\psi'(w)}{\psi(w)-z}=\sum_{n=0}^\infty
\frac{F_n(z)}{w^{n+1}}.
\end{equation}

By an application of Cauchy integral formula,
(\ref{FabEq84}) yields the following integral
representation for the Faber polynomials: for
every $R>1$ and $z$ lying in the interior of the
level curve $L_R:=\{\psi(w):|w|=R$\},
\begin{equation}\label{FabEq27}
F_n(z)=\frac{1}{2\pi
i}\oint_{\mathbb{T}_R}\frac{t^n\psi'(t)dt}{\psi(t)-z}\,,
\end{equation}
while for $z$ lying in the exterior of $L_R$,
\begin{equation}\label{FabEq78}
F_n(z)=[\phi(z)]^n+\frac{1}{2\pi
i}\oint_{\mathbb{T}_R}\frac{t^n\psi'(t)dt}{\psi(t)-z}.
\end{equation}

It this paper we investigate the asymptotic behavior of the
Faber polynomials and their zeros for certain domains
$\Omega$ that are bounded by piecewise analytic curves.
Hence the title of this paper. More precisely, we consider
domains (or equivalently, functions $\psi$) satisfying
assumptions A.1 and A.2 to be stated in what follows.

\begin{figure}
\centering
\includegraphics[scale=.6]{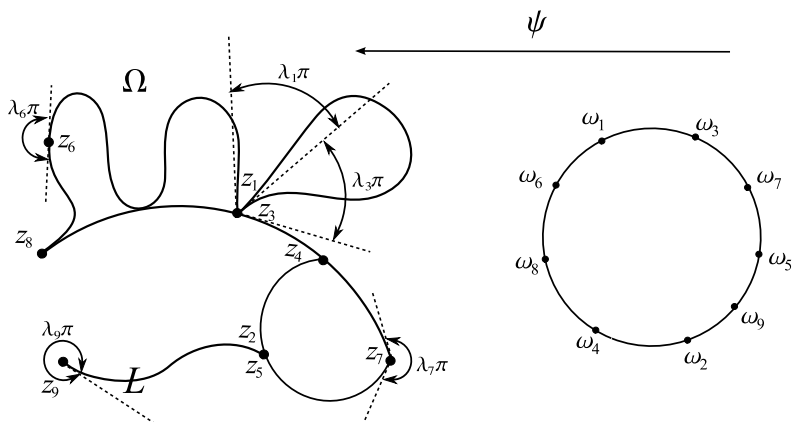}
\caption{Illustration of a map $\psi$ satisfying conditions
A.1 and A.2.}\label{FabFig4}
\end{figure}

We define an \emph{analytic arc} as being the image of the
interval $[0,1]$ by a function $f(t)$ analytic in $[0,1]$
such that $f'(t)\not=0$ for all $t\in [0,1]$ and
$f(t_1)\not=f(t_2)$ for all $0<t_1<t_2< 1$. The endpoints
of the arc are $f(0)$ and $f(1)$, which may coincide. We
call the arc \emph{simple} if $f$ is one-to-one on $[0,1]$.
Notice that, according to this definition, an analytic
Jordan curve is also an analytic arc. Our first assumption
is:

\begin{enumerate}
\item[\textbf{A.1:}] The map $\psi$ has a continuous extension to $\cj{\Delta}_1$
and there are $s\geq 1$ distinct points
$\ou_1,\ou_2,\ldots,\ou_s$ in $\mathbb{T}_1$ such that if
$\ell$ is any of the $s$ open circular arcs that compose
$\mathbb{T}_1\setminus \{\ou_1,\ou_2,\ldots,\ou_s\}$, say
with endpoints $\ou_k$, $\ou_j$, then $\psi$ is one-to-one
on $\ell$ and $\psi\left(\cj{\ell}\right)$ is an analytic
arc with endpoints $\psi(\ou_k)$, $\psi(\ou_j)$ (see Figure
\ref{FabFig4}).
\end{enumerate}

Thus, $\partial \Omega$ is a  piecewise analytic
curve that we denote by $L$. Let $z\in L$ and
$w=e^{i\Theta}$ be such that $z=\psi(w)$. The
\emph{exterior angle at $z$ relative to $w$} is
defined to be that angle $\alpha\in [0, 2\pi]$
such that
\[
\arg\left[\psi\left(e^{i\theta}\right)-z\right]\to\left\{\begin{array}{ll}
                                   \beta &\ \,\mathrm{as}\ \theta\to \Theta-\,,\\
                                   \beta+\alpha &\ \,\mathrm{as}\
                                   \theta\to\Theta+\,.
                                 \end{array}
\right.
\]

Let
\[
z_k:=\psi(\ou_k), \quad k\in\{1,2,\ldots,s\}.
 \]
These points $z_k$ will be called the corners of $L$.
Notice that they are not necessarily pairwise distinct.

For each $ k\in\{1,2,\ldots,s\}$, let
$\lambda_k\in [0,2]$ be such that $\lambda_k\pi$
is the exterior angle at $z_k$ relative to
$\ou_k$. It is well-known that when
$\lambda_k>0$, the mapping $\psi$ has an
asymptotic expansion about $\ou_k$ in functions
of the form
\begin{equation}\label{FabEq67}
(w-\ou_k)^{l+j\lambda_k}(\log(w-\ou_k))^m\,,\quad w\in
\mathbb{C}\setminus\{t\ou_k:t\leq 1\},
\end{equation}
with $l\geq 0$, $j\geq1$, and $m\geq 0$ integers (see
\cite{Lehman} and also \cite[pp. 57-58]{Pomm}). We will
refer to it as the Lehman expansion of $\psi$ about $\ou_k$
and its exact meaning is explained in Section \ref{FabSec1}
below. Logarithmic terms (i.e., functions of the form
(\ref{FabEq67}) with $m\geq 1$) may occur in the expansion
only if $\lambda_k$ is a rational number. Our second
assumption on $\psi$ is:
\begin{enumerate}
\item[\textbf{A.2:}] $\lambda_k>0$ for every $ k\in\{1,2,\ldots,s\}$, and if
$\lambda_k\in\{1,2\}$ for all $
k\in\{1,2,\ldots,s\}$, then there is at least one
$k$ for which logarithmic terms occur in the
Lehman expansion of $\psi$ about $\ou_k$.
\end{enumerate}

If $\lambda_k\not\in \{1,2\}$, then $\ou_k$ is a
singularity of $\psi$, and a sufficient condition for an
$\ou_k$ with $\lambda_k\in\{1,2\}$ to be a singularity of
$\psi$ is precisely that logarithmic terms occur in the
Lehman expansion of $\psi$ about $\ou_k$. We do not know
whether this condition is also necessary. If that were the
case, we could simply rephrase A.2 by saying that all the
$\lambda_k$'s are positive and at least one $\ou_k$ is a
singularity of $\psi$.

Let us then consider a map $\psi$ satisfying A.1 and A.2.
The letter $G$ will denote the complement of $\cj{\Omega}$,
so that if $L$ is a Jordan curve, $G$ is the interior
domain of $L$.

A first observation is that the asymptotic
behavior of $F_n$ in $\Omega$ is already given by
the integral representation (\ref{FabEq78}): for
arbitrary $1<R<r$,
\begin{equation}\label{FabEq88}
F_n(z)=[\phi(z)]^n\left[1+\mathcal{O}\left(\frac{r^n}{R^n}\right)\right]\,,
\end{equation}
uniformly on $\psi\left(\cj{\Delta}_r\right)$ as
$n\to\infty$. Hence, every closed subset of $\Omega$ will
be free of zeros of $F_n$ for $n$ large enough, and all
accumulation points\footnote{$t$ is an accumulation point
if every neighborhood of $t$ contains zeros of infinitely
many polynomials $F_n$.} of the zeros of the Faber
polynomials must be contained in $\cj{G}$.

Formula (\ref{FabEq88}) has been previously extended to $L$
in the pointwise sense, under the additional assumption
that $L$ is a Jordan curve. In this case $\phi$ has a
continuous extension to $L$ and a more general result of
Pritsker \cite[Thm. 1.1]{Pritsker1} about the behavior of
weighted Faber polynomials implies that if $z\in L$ is not
a corner, then
\begin{equation}\label{FabEq60}
F_n(z)=[\phi(z)]^n(1+o(1)),\quad (n\to\infty)
\end{equation}
while for every corner $z_k$,
\begin{equation}\label{FabEq50}
F_n(z_k)=\lambda_k[\phi(z_k)]^n(1+o(1)),\quad (n\to\infty).
\end{equation}

The behavior of $F_n$ in $G$ has remained quite
unknown, but at least for $L$ a piecewise
analytic Jordan curve without cusps, Gaier
\cite{Gaier} was able to derive uniform estimates
on the decrease of $F_n$ of the form
\begin{equation} \label{FabEq102}
F_n(z)=\mathcal{O}\left(n^{-\lambda}\right),\quad z\in G,
\end{equation}
where $\lambda$ is the smallest of the exterior angles at
the corners of $L$.

In this paper we much improve these results by
providing asymptotic formulas for $F_n$  that do
not require $L$ to be a Jordan curve and that
hold \emph{uniformly} on closed subsets of the
complex plane. Moreover, our estimates for the
rate of decay of the error terms involved are, in
general, best possible. Theorems \ref{FabThm2},
\ref{FabThm8} and \ref{FabThm9} of Section
\ref{FabSubSec2} are the strengthened versions of
(\ref{FabEq88}), (\ref{FabEq60}) and
(\ref{FabEq50}), while Theorem \ref{FabThm1}
transparently describes the behavior of $F_n$ in
$G$, yielding, in particular, Gaier's estimate
(\ref{FabEq102}). Theorem \ref{FabThm1} also
shows that the pointwise estimates given by Gaier
in \cite[Thm. 2]{Gaier}
 indeed hold locally uniformly on $G$.

As for the zeros of $F_n$, the fact that the map
$\psi$ under consideration has a singularity on
$\mathbb{T}_1$ implies, by a general result of
Ullman \cite[Thm. 1]{Ullman}, that all points of
$L$ are accumulation points of the zeros of the
Faber polynomials. A later complement to Ullman's
results by Kuijlaars and Saff \cite[Thms. 1.3,
1.4]{KS} implies that there is always a
subsequence of the sequence $\{\nu_{n}\}_{n\geq
1}$  of normalized counting measures of the zeros
of the $F_n$'s that converges in the
weak*-topology to the equilibrium measure $\mu_L$
of $L$, and this is true of the entire sequence
provided that $G=\emptyset$ (see (\ref{FabEq28})
and (\ref{FabEq77}) in Section \ref{FabSubSec1}
for definitions of $\nu_n$ and $\mu_L$).

We will be able to say much more. In Section
\ref{FabSubSec1} we show that, independently of
whether $G$ is connected or not, there is always
a subsequence of $\{\nu_{n}\}_{n\geq 1}$ that
converges in the weak*-sense to $\mu_L$. In fact,
under an additional assumption that is naturally
satisfied in a large number of cases (including
when $L$ is a Jordan curve), we prove that
compact subsets of $G$ contain at most a finite
(independent of $n$) number of zeros of every
$F_n$, forcing the whole sequence
$\{\nu_{n}\}_{n\geq 1}$ to converge to $\mu_L$.
Furthermore, under that assumption we are also
able to characterize those points of $G$ that are
accumulation points of the zeros of Faber
polynomials.

Faber polynomials for particular domains of the
complex plane has been the subject of many recent
works, in several of which the boundary of the
domain is precisely a piecewise analytic curve,
for example, $m$-stars \cite{He3}, \cite{K},
circular lunes \cite{He1}, $m$-fold symmetric
curves and certain lemniscates \cite{He2},
annular and circular sectors \cite{Coleman1},
\cite{Coleman2}. Our results apply to all these
examples.

The rest of this paper is organized as follows. Section
\ref{FabSubSec2} presents the asymptotic formulas for Faber
polynomials, in Section \ref{FabSubSec1} we draw some
conclusions on their zero behavior and analyze two concrete
examples. In Section \ref{FabSec1} we discuss in detail the
Lehman expansion of the exterior map $\psi$, and finally in
Sections \ref{FabSec2} and \ref{FabSec3} we prove all the
results.

\section{Asymptotic behavior of
$F_n(z)$}\label{FabSubSec2}

Recall that we are considering a map $\psi$
satisfying assumptions A.1 and A.2 stated in the
introduction. Let
$\Theta_1,\Theta_2,\ldots,\Theta_s$ be the
arguments of the numbers $\ou_k$, that is,
\[
\ou_k=e^{i\Theta_k},\quad 0\leq\Theta_k<2\pi,\quad 1\leq
k\leq s\,.
\]

Assumption A.2 is independent of the branches
chosen for the functions in (\ref{FabEq67}) in a
$\delta$-neighborhood of the form
$\{w\in\Delta_1:0<|w-\ou_k|<\delta\}$. However,
to simplify the statements of our results, we
choose those corresponding to the branch of the
argument
\[
\Theta_k-\pi<\arg(w-\ou_k)<\Theta_k+\pi,\quad
w\in \mathbb{C}\setminus\{t\ou_k:t\leq 1\}.
\]

We shall say that $\ou_k$ is \emph{relevant} if
either $\lambda_k\not\in\{0,1,2\}$, or  if
$\lambda_k\in\{1,2\}$ and logarithmic terms occur
in the expansion of $\psi$ about $\ou_k$. With
this definition, condition A.2 states that all
$\lambda_k$'s are positive and that \emph{there
is at least one relevant $\ou_k$}.

Let now $v\geq 1$ be the number of relevant
$\ou_k$'s. Hereafter we shall assume that the
$\ou_k$'s have been indexed in such a way that
$\ou_1,\ou_2,\ldots, \ou_v$ are precisely the
relevant ones. The following weaker version of
the Lehman expansion of $\psi$ about a relevant
$\ou_k$ is sufficient to state our main results.

If $k\in\{1,2,\ldots,v\}$ is such that
$\lambda_k\not\in\{1,2\}$, then there is
$A_k\not=0$ such that as $w\to \ou_k$ from the
exterior of the unit circle,
\begin{equation}\label{FabEq3}
\psi(w)=z_k+A_k (w-\ou_k)^{\lambda_k}\left(1+o(1)\right),
\end{equation}
while if $k\in\{1,2,\ldots,v\}$ is such that
$\lambda_k\in\{1,2\}$, then there exist positive
integers $r_k,\, m_k$ with $r_k\geq \lambda_k$
and  $1\leq m_k\leq \lfloor
r_k/\lambda_k\rfloor$, and complex numbers
$A_k\not=0$, $c^k_0\not=0$,
$c^k_1,c^k_2,\ldots,c^k_{r_k-1}$, such that as
$w\to \ou_k$ from the exterior of the unit
circle,
\begin{eqnarray}\label{FabEq32}
\psi(w)&=&z_k+\sum_{l=0}^{r_k-1}c^{k}_{l}
(w-\ou_k)^{\lambda_k+l}+A_k(w-\ou_k)^{r_k+\lambda_k}\left(\log
(w-\ou_k)\right)^{m_k}\left(1+o(1)\right).
\end{eqnarray}

From relations (\ref{FabEq3}) and (\ref{FabEq32}), we associate to
each relevant $\ou_k$ ($1\leq k\leq v$) the number $A_k$, the
numbers $r_k$ and $m_k$ whenever $\lambda_k\in\{1,2\}$, and the
following pair:
\begin{eqnarray*}
(\Lambda_k,M_k):=\left\{\begin{array}{cc}
                          (\lambda_k,0), &\ \mathrm{if}\ \lambda_k\not\in \{1,2\}, \\
                          (r_k+\lambda_k,m_k-1), &\ \mathrm{if}\ \lambda_k\in
                          \{1,2\}.
                        \end{array}
\right.
\end{eqnarray*}
Observe that $\Lambda_k\geq 2$ if $\lambda_k\in\{1,2\}$.  We will
say that $(\Lambda_k,M_k)<(\Lambda_j,M_j)$ if either
$\Lambda_k<\Lambda_j$, or $\Lambda_k=\Lambda_j$ and $M_k>M_j$.

By reindexing the relevant $\ou_k$'s if needed,
we may assume that $\ou_1,\ldots,\ou_v$ are such
that
\[
\left(\Lambda_1,M_1\right)=\cdots=\left(\Lambda_u,M_u\right)<\left(\Lambda_{u+1},M_{u+1}\right)\leq
\cdots \leq \left(\Lambda_v,M_v\right),
\]
for some $u\in\{1,2,\ldots,v\}$.

For any two integers $n,\,m\geq 0$ and a real
$\beta>-1$, we define
\begin{equation}
\alpha_{\beta,m}(n):=\int_0^1x^n(1-x)^\beta\left(\log
(1-x)\right)^m dx\,.
\end{equation}
Then,\footnote{ A proof of (\ref{FabEq16}) is
given at the end of Section \ref{FabSec2}.}
\begin{eqnarray}\label{FabEq16}
\alpha_{\beta,m}(n)&=&\frac{\Gamma(\beta+1)n!(-\log
n)^m\left[1+\mathcal{O}\left(1/\log n\right)\right]}{\Gamma(n+\beta+2)}\\
&=& \frac{\Gamma(\beta+1)(-\log
n)^m(1+o(1))}{n^{\beta+1}}\nonumber
\end{eqnarray}
as $n\to\infty$.

Recall that we have defined $L:=\partial\Omega$
and $G:=\mathbb{C}\setminus \cj{\Omega}$. We
first consider the behavior of $F_n$ on $G$.

\begin{thm}\label{FabThm1} Let\,\footnote{The letter $\Gamma$ stands for the Euler gamma function.}
\[
\mathcal{C}_1:=\left\{\begin{array}{ll}
                        [\Gamma(-\lambda_1)\Gamma(\lambda_1)]^{-1}, &\quad if\ \lambda_1\not\in\{1,2\}, \\
                        (-1)^{\Lambda_1-1}m_1\Lambda_1, & \quad if\
                        \lambda_1\in\{1,2\}.
                      \end{array}
\right.
\]
Then, for every $z\in G$,
\begin{equation}\label{FabEq22}
\frac{F_n(z)}{\alpha_{\Lambda_1-1,M_1}(n)}=\mathcal{C}_1\sum_{k=1}^u
\frac{A_k e^{i(n+\Lambda_1)\Theta_k}}{z-z_k}+R_n(z),
\end{equation}
where $R_n(z)$ converges to zero locally
uniformly  on $G$.
\end{thm}

\begin{rem}\label{FabRem1}\emph{Since some of the $z_k$'s may coincide, it is possible that
for some subsequence $\{n_j\}\subset\mathbb{N}$,
the rational functions $\sum_{k=1}^u A_k
e^{i(n_j+\Lambda_1)\Theta_k}/(z-z_k)$ occurring
in (\ref{FabEq22}) be (or at least converge to)
the constant zero function (see the example
discussed at the end of Section
\ref{FabSubSec1}). Nevertheless, as we show with
that example (see Theorem \ref{FabThm3} and its
proof), in a situation like this we could still
be successful in proving that, after proper
normalization, $\{F_{n_j}\}$ behaves like certain
sequence of rational functions that do not
approach zero. The proof can be attempted as
follows: write (\ref{FabEq27}) as in
(\ref{FabEq34}), then combine identity
(\ref{FabEq5}) with the Lehman expansion of
$\psi$ about $\ou_k$ to obtain, for each of the
integrals under the $\Sigma$ sign of
(\ref{FabEq34}), subsequent terms of its
expansion as a sum of rational functions whose
denominators are powers of $(z-z_k)$. }
\end{rem}

Let us now turn our attention to $\cj{\Omega}$.
For every $z\in \cj{\Omega}$, let $\eta(z)$ be
the (finite) number of elements of the set
$\left\{w\in \cj{\Delta}_1:\psi(w)=z\right\}$.
These elements will be denoted by
\[\phi_1(z),\phi_2(z),\ldots,\phi_{\eta(z)}(z),\]
being irrelevant the order in which they are
numerated. Of course, $\eta(z)=1$ and
$\phi_1(z)=\phi(z)$ for $z\in \Omega$.

Because $\Omega$ has no outward-pointing cusps,
if $z\in L$ and two elements of $\{w\in
\mathbb{T}_1:\psi(w)=z\}$ belong to
$\mathbb{T}_1\setminus\{\ou_1,\ldots,\ou_s\}$,
then indeed $\eta(z)=2$. Hence, $1\leq
\eta(z)\leq s+1$, and when $z$ is not a corner,
$\eta(z)\leq 2$.

For every $z\in L$ and $1\leq j\leq \eta(z)$, let
$\hat{\lambda}_j(z)\pi$ ($0< \hat{\lambda}_j(z)\leq 2$) be
the exterior angle at $z$ relative to $\phi_j(z)$. Then,
only if $z$ is a corner of $L$ it is possible to have
$\hat{\lambda}_j(z)\not=1$ for some $1\leq j\leq \eta(z)$.
Let us define
\begin{equation}\label{FabEq69} L_1:=\left\{z\in L: \eta(z)=1,\ \hat{\lambda}_1(z)=1\right\}
,
\end{equation}
\begin{equation}\label{FabEq91}
L_2:=\left\{z\in L: \eta(z)=2,\
\hat{\lambda}_1(z)=\hat{\lambda}_2(z)=1
\right\}\cup\left\{z_k: \lambda_k=2\right\}\,.
\end{equation}

Observe that \[L=L_1\cup L_2\cup
\{z_1,\ldots,z_s\}.
\]

Let $\Phi_n:\cj{\Omega}\to\cj{\mathbb{C}}$ be
defined by
\begin{equation*}
    \Phi_n(z):=\left\{\begin{array}{ll}
                        \underset{}{[\phi(z)]^n}\,, &\ \,\ z\in\Omega,
                         \\
                       {\D \sum_{j=1}^{\eta(z)}\hat{\lambda}_j(z)[\phi_j(z)]^n}\,, &
                        \ \,\ z\in L.
                      \end{array}
    \right.
\end{equation*}
\begin{thm}\label{FabThm2} For every $z\in \cj{\Omega}$,
\begin{equation}
F_n(z)=\Phi_n(z)+R_n(z),
\end{equation}
where $R_n(z)$ is such that
\begin{enumerate}
\item[(a)] if $E$ is a closed set and either $E\subset\Omega\cup
L_1$ or $E\subset L_2$, then
$R_n(z)=\mathcal{O}\left(n^{-\Lambda^*}(\log
n)^{M^*}\right)$ uniformly as $n\to\infty$ on
$E$, where $(\Lambda^*,M^*)$ is the smallest
element of the set
\[\left\{(\Lambda_1,M_1)\right\}\cup\left\{(r_k,M_k):1\leq k\leq
u,\ z_k\in E\right\};\]

\item[(b)] for $j=1,2,\ldots,s$,
$R_n(z_j)=\mathcal{O}\left(n^{-\Lambda^*_j}(\log
n)^{M^*_j}\right)$ as $n\to\infty$, where
$(\Lambda^*_j,M^*_j)$ is the smallest element of
the set
\[\left\{(\Lambda_1,M_1)\right\}\cup\left\{(r_k,M_k):\psi(\ou_k)=z_j,\
\lambda_k\in \{1,2\}\right\}.\]
\end{enumerate}
\end{thm}

We can be more specific for closed subsets of
$\Omega\cup L_1$ or $L_2$ without corners.
\begin{thm}\label{FabThm8} Let $E\subset \left(\Omega\cup L_1\right)\setminus \{z_1,\ldots,
z_s\}$ be a closed set. There exists an open set
$U\supset E$ such that $\phi$ has an analytic and
univalent continuation to $U$ and
\begin{equation}\label{FabEq76}
F_n(z)=[\phi(z)]^n+\alpha_{\Lambda_1-1,M_1}(n)\left(\mathcal{C}_1\sum_{k=1}^u
\frac{A_k
e^{i(n+\Lambda_1)\Theta_k}}{z-z_k}+R_n(z)\right),
\end{equation}
with $R_n(z)\to 0$ uniformly on $U$ as $n\to\infty$.
\end{thm}

Let now $E\subset L_2\setminus\{z_1,\ldots,z_s\}$
be a simple analytic arc. Then, there is a
``strip-like" connected neighborhood $U$ of $E$
such that $\cj{U}\cap L$ is a simple analytic arc
contained in $L_2$, $U\setminus L$ consists of
two open components $U^+$, $U^-$, both contained
in $\Omega$, and if $z^*$ denotes the Schwarz
reflection of $z$ about the analytic arc $E$,
then $z^*\in U^{\pm}$ if and only if $z\in
U^{\mp}$.

Let $\phi_{+}$, $\phi_{-}$ be the restrictions of
$\phi$ to $U^+$, $U^-$, respectively. Each of
these functions is continuous along the arc
$U\cap L$, mapping it onto an arc of the unit
circle. By the Schwarz reflection principle
\cite{Davis}, each function has an analytic and
univalent continuation to all of $U$, whose
values on $U^{\pm}$ are given by
\[
\phi_{\pm}(z)=\frac{1}{\cj{\phi_{\pm}(z^*)}},\quad z\in
U^{\mp}.
\]

\begin{thm}\label{FabThm9}Let $E\subset L_2\setminus\{z_1,\ldots,z_s\}$ be a simple analytic arc. There
exists a neighborhood $U$ of $E$ as described
above such that for all $z\in U$,
\begin{equation}\label{FabEq87}
F_n(z)=[\phi_+(z)]^n+[\phi_-(z)]^n+\alpha_{\Lambda_1-1,M_1}(n)\left(\mathcal{C}_1\sum_{k=1}^u
\frac{A_k e^{i(n+\Lambda_1)\Theta_k}}{z-z_k}+R_n(z)\right),
\end{equation}
with $R_n(z)\to 0$ uniformly on $U$ as $n\to\infty$.
\end{thm}

\begin{rem}\label{FabRem3}\emph{1) Concerning how fast the error terms $R_n(z)$
in (\ref{FabEq22}), (\ref{FabEq76}) and
(\ref{FabEq87}) approach zero, the best it can be
said, in general, is that they decrease at least
as fast as the dominant terms in the right-hand
side of (\ref{FabEq82}) in page
\pageref{FabEq82}, where the rate of decay of the
functions $r_{\sigma_k,n}(z)$ therein is
estimated in the table of Remark \ref{FabRem2} in
page \pageref{FabPag1}.}

\emph{2) The estimates provided in Theorem
\ref{FabThm2} for $R_n(z)$ are also best
possible, as can be verified from relation
(\ref{FabEq24}) for part (a), and from relation
(\ref{FabEq55}) for part (b).}

\end{rem}

\section{The zeros of
$F_n(z)$}\label{FabSubSec1}

In this section we draw from our previous results
some conclusions about the location, accumulation
points and limiting distribution of the zeros of
Faber polynomials.

From Theorem \ref{FabThm2} we immediately see
that
\begin{equation}\label{FabEq68}
\lim_{n\to\infty}\frac{F_n(z)}{\Phi_n(z)}=1
\end{equation}
locally uniformly on $\Omega\cup L_1$, where
$\Phi_n|_{\Omega\cup L_1}$ is simply the
continuous extension of $[\phi(z)]^n$ to
$\Omega\cup L_1$. Hence, we have

\begin{cor}\label{FabCor2}For every closed set $E\subset \Omega\cup
L_1$, there is a number $N_E$ such that when
$n>N_E$, $F_n(z)$ has no zeros on $E$.
\end{cor}

Let us now focus on the effect that Theorem
\ref{FabThm1} has on the zeros of $F_n$. It is
interesting that asymptotic formulas similar to
(\ref{FabEq22}) are also satisfied by orthogonal
polynomials on the unit circle with respect to
certain types of weights. Some of the results
that follow are basically known consequences of
such type of behavior, see e.g., \cite{Sza},
\cite{andrei}.

We first rewrite (\ref{FabEq22}) in a more
suitable way. Put
\[
\hat{A}_k:=A_ke^{i\Lambda_1(\Theta_k-\Theta_1)}\,,\quad
1\leq k\leq u\,,
\]
and let $\theta_1,\theta_2,\ldots,\theta_u$ be such that
\[
e^{2\pi i\theta_k}=e^{i(\Theta_k-\Theta_1)},
\quad \theta_k\in (0,1],\quad 1\leq k\leq u\,,
\]
so that (\ref{FabEq22}) takes the form
\begin{equation}\label{FabEq64}
F^*_n(z)=H_n(z)+o(1)
\end{equation}
locally uniformly on $G$ as $n\to\infty$, where
\begin{equation}\label{FabEq85}
F_n^*(z)=\left[\mathcal{C}_1e^{i
(n+\Lambda_1)\Theta_1}\alpha_{\Lambda_1-1,M_1}(n)\right]^{-1}F_n(z)\,,\quad
H_n(z)=\sum_{k=1}^u \frac{\hat{A}_k e^{2\pi
in\theta_k}}{z-z_k}\,.
\end{equation}

In view of (\ref{FabEq64}) and the form of the
rational functions $H_n$, the sequence $\{F_n^*\}
_{n\geq 1}$ is a normal family on $G$, and a
function $f$ is the uniform limit on $G$ of some
subsequence $\left\{F_{n_j}^*\right\} _{j\geq 1}$
if and only if it is the uniform limit of
$\left\{H_{n_j}\right\} _{j\geq 1}$. Hence, every
such $f$ must have the form
\begin{equation}\label{FabEq81}
f(z)=\sum_{k=1}^u \frac{\hat{A}_k e^{2\pi
i\vartheta_k}}{z-z_k}.
\end{equation}

Because the $z_k$'s are not necessarily pairwise distinct,
some of these limit functions can be identically zero,
which makes Theorem \ref{FabThm1} insufficient to describe
the zeros of the $F_n$'s. Therefore, we shall often make
the assumption that
\begin{enumerate}
\item[\textbf{A.3:}] no subsequence of $\{H_n\}_{n\geq 0}$ converges
to the null function.
\end{enumerate}

If A.3 is satisfied, then all uniform limit
points of $\{F^*_n\}_{n\geq 0}$ are nonzero
rational functions of bounded degree. Let us see
what this implies on the limiting distribution of
the zeros of $F_n$.

Let $\nu_{n}$ be the normalized counting measure
of the zeros of $F_n$, that is,
\begin{equation}\label{FabEq28}
\nu_n:=n^{-1}\sum_{k=1}^n\delta_{z_{k,n}}\,,\quad
n=1,2,\ldots,
\end{equation}
where $z_{1,n},z_{2,n},\ldots,z_{n,n}$ are the
zeros of $F_n$ (counting multiplicities) and
$\delta_z$ is the unit point measure at $z$.

A subsequence $\{\nu_{n_j}\}_{j\geq 1}$ of
$\{\nu_{n}\}_{n\geq 1}$ is said to converge in
the weak*-topology to a Borel measure $\mu$
(symbolically, $\nu_{n_j}\wc \mu$ as
$j\to\infty$) if for every continuous function
$f$ defined on $\cj{\mathbb{C}}$,
$\lim_{j\to\infty}\int f\nu_{n_j}=\int f\mu$.

Let $\mu_L$ be the \emph{equilibrium measure} of
$L$, i.e., the measure supported on $L$ whose
value at any given Borel set $B\subset L$ is
\begin{equation}\label{FabEq77}
\mu_L(B)=\frac{1}{2\pi}\int_{B^{-1}}|dt|,\quad
B^{-1}:=\{t\in \mathbb{T}_1: \psi(t)\in B\}.
\end{equation}
Notice that $\mu_L$ is a probability measure
whose support is $L$.

\begin{cor}\label{FabPro1} Assume that A.3 holds. Then, for every closed set $E\subset G$ there is a
number $N_E$ such that when $n>N_E$, $F_n(z)$ has
at most $J-1$ zeros in $E$ (counting
multiplicities), where $J$ is the number of
corners $z_k$. Hence, $\nu_{n}\wc \mu_L$ as
$n\to\infty$.
\end{cor}

\begin{rem}\emph{Under assumption A.3, finer
results similar to Thm. 4 of \cite{andrei} (see
also \cite[Thms. 11.1, 11.2]{Simon}) on the
separation, distribution and speed of convergence
to $L$ of those zeros of $F_n$ that lie near
$L_1$ can be derived from Theorem \ref{FabThm8}.
}
\end{rem}

Condition A.3 holds in a large number of cases. For
instance, if there is $k$ such that $z_j\not= z_k$ whenever
$j\not =k$, as is the case of $L$ a Jordan curve. Indeed,
if A.3 does not hold, there must be a limit function $f$
(which has the form (\ref{FabEq81})) such that for all
$1\leq j\leq u$, $\sum_{k\,:\,z_k=z_j}\hat{A}_ke^{2\pi
i\vartheta_k}=0$. Certain numbers $\vartheta_k$ satisfying
this last equality can be found if and only if
\[
2\max_{k\,:\,z_k=z_j}|A_{k}|\leq
\sum_{k\,:\,z_k=z_j}|A_{k}|.
\]
However, whether these found $\vartheta_k$'s
actually correspond to a limit function $f$
depends on the specific values of the
$\vartheta_k$'s and can be determined from the
general form of the uniform limit points of
$\{H_n\}_{n\geq 0}$ that we establish next.

Among the numbers
$1=\theta_1,\theta_2,\ldots,\theta_u$, there is a
basis over the rationals  containing $\theta_1$
\cite[Ch. III. p. 4]{Cassels}, say $\theta_1,
\theta_2,\ldots,\theta_{u^*}$, $1\leq u^*\leq u$,
such that for every $k\in \{1,2,\ldots,u\}$,
there are unique rational numbers
$r_{k1},r_{k2},\ldots,r_{ku^*}$ with
\[
\theta_k=\sum_{j=1}^{u^*}r_{kj}\theta_j, \quad 1\leq k \leq
u.
\]
Notice that $u^*=1$ if and only if all the $\theta_k$'s are
rational, and if $u^*\geq 2$, then
$\theta_2,\ldots,\theta_{u^*}$ are irrational numbers
linearly independent over the rationals.

For every $k\in \{1,2,\ldots,u\}$, let $1\leq p_k
\leq q_k$ be the unique relatively prime integers
such that
\[
e^{2\pi i\,r_{k1}}=e^{2\pi i\,p_k/q_k},
\]
so that
\begin{equation}\label{FabEq99}
e^{2\pi i\theta_k}=e^{2\pi
i\,\left(\frac{p_k}{q_k}+\sum_{j=2}^{u^*}r_{kj}\theta_j\right)},
\quad 1\leq k \leq u,
\end{equation}
where in case $u^*=1$, the sum
$\sum_{j=2}^{u^*}\cdots$ above is understood to
be zero (notice that $p_1=q_1=1$, but $p_k<q_k$
for $k>1$).

Let $\mathrm{\mathbf{q}}$ be the least common
multiple of the denominators
$q_1,q_2,\ldots,q_u$, and for every
$\ell=\{1,2,\ldots, \mathrm{\mathbf{q}}\}$, let
\[
\ell p_k=s_{k\ell}\!\!\!\mod q_k,\quad 0\leq s_{k\ell}<
q_k\,.
\]
Observe that  two $u$-tuples
$\left(s_{1\ell},s_{2\ell},\ldots,s_{u\ell}\right)$
corresponding to different values of $\ell$ are distinct.

\begin{thm} \label{FabThm4} The functions $f$ that are the uniform limit of some subsequence of $\{H_n\}_{n\geq 0}$
are the rational functions of the form
\begin{equation}\label{FabEq100}
f(z)=\sum_{k=1}^u \frac{\hat{A}_k e^{2\pi
i\,\left(\frac{s_{k\ell}}{q_k}+\sum_{j=2}^{u^*}r_{kj}\alpha_j\right)}}{z-z_k}
\end{equation}
with $\ell=\{1,2,\ldots, \mathrm{\mathbf{q}}\}$
and $\alpha_2,\ldots,\alpha_{u^*}$ arbitrary real
numbers. In particular, there is always such a
limit function $f$ that is not identically zero.
\end{thm}

As mentioned in the introduction, a result of
Kuijlaars and Saff \cite[Thms. 1.3, 1.4]{KS}
implies that if $G$ is connected, then some
subsequence of the counting measures
$\{\nu_{n}\}$ must converge in the weak*-sense to
the equilibrium measure of $L$. From Theorem
\ref{FabThm4} we now see that the connectedness
of $G$  can be dropped.

\begin{cor}\label{FabCor1}
There is always a subsequence $\{n_j\}\subset
\mathbb{N}$ such that $\nu_{n_j}\wc \mu_L$ as
$j\to\infty$.
\end{cor}

In fact, we have seen that as long as A.3 is satisfied
(even if $G$ is disconnected), it is true that $\nu_{n}\wc
\mu_L$ as $n\to\infty$. But there are examples with $G$
disconnected and some subsequence of $\{\nu_{n}\}$
converging to a measure supported in $G$ (see the example
discussed at the end of this section). We have not been
able to determine, however, whether the connectedness of
$G$ is sufficient for $\nu_{n}\wc \mu_L$ as $n\to\infty$.
We leave it as a

\begin{conj} If $G$ is connected, then
$\nu_{n}\wc \mu_L$ as $n\to\infty$.
\end{conj}

Let us now concentrate on the set $\mathcal{Z}$
of accumulation points of the zeros of the Faber
polynomials, i.e., $\mathcal{Z}$ is the set of
all points $t\in \cj{\mathbb{C}}$ such that every
neighborhood of $t$ contains zeros of infinitely
many polynomials $F_n$.

As we pointed out in the introduction, it is
always the case that $\Omega\cap
\mathcal{Z}=\emptyset$, and having the maps
$\psi$ under consideration a singularity on
$\mathbb{T}_1$, a general result of Ullman
\cite[Thm. 1]{Ullman} implies that $L\subset
\mathcal{Z}$ (this also follows from Corollary
\ref{FabCor1} since the support of $\mu_L$ is
$L$). The following characterization of
$\mathcal{Z}\cap G$ follows directly from Theorem
\ref{FabThm4} and Hurwitz's Theorem.

\begin{cor}\label{FabCor3} Assume A.3
holds. The point $t\in G$ also belongs to
$\mathcal{Z}$ if and only if there exist an
integer $\ell=\{1,2,\ldots,
\mathrm{\mathbf{q}}\}$ and real numbers
$\alpha_2,\ldots,\alpha_{u^*}$  such that
\begin{equation}\label{FabEq59}
\sum_{k=1}^u \frac{\hat{A}_k e^{2\pi
i\,\left(\frac{s_{k\ell}}{q_k}+\sum_{j=2}^{u^*}r_{kj}\alpha_j\right)}}{t-z_k}=0.
\end{equation}
\end{cor}

\begin{rem} \emph{Assume A.3 holds, so that by Corollary \ref{FabCor3} we have the following. If $z_1=z_2=\cdots=z_u$, then $\mathcal{Z}\cap G=\emptyset$. Otherwise:
\begin{enumerate}
\item[a)] if $u^*=1$ (i.e., all the $\theta_j$'s are
rational), then the number of points in $\mathcal{Z}\cap G$
is finite, namely at most $(u-1)\mathrm{\mathbf{q}}$;
\item[b)] if $u^*=2$, then by fixing $\ell$ and letting $\alpha_2$ vary,
equation (\ref{FabEq59}) can be written as
\begin{equation}\label{FabEq80}
g_{0,\ell}(z)+g_{1,\ell}(z)t+\cdots
+g_{u-1,\ell}(z)t^{u-1}=0, \quad |z|=1,
\end{equation}
where the $g_{k,\ell}(z)$'s are certain
polynomials, so that if $f_1,\ldots,f_m$ are the
algebraic functions determined by the algebraic
equations in (\ref{FabEq80}) (see e.g.,
\cite[Chap. 5]{Knopp}), then $\mathcal{Z}\cap G$
consists of the traces left in $G$ by the curves
$f_1(\mathbb{T}_1),\ldots,f_m(\mathbb{T}_1)$,
plus possibly some of the solution points
corresponding to the algebraic singularities of
the $f_k$'s. In particular, when $u=2$, equation
(\ref{FabEq59}) reduces to $
|\hat{A}_1(t-z_2)|=|\hat{A}_2(t-z_1)|$, so that
$\mathcal{Z}\cap G$ is the trace in $G$ of a line
if $|A_1|=|A_2|$,  or of a circle if
$|A_1|\not=|A_2|$;
\item[c)] if $u^*\geq 2$, then $\mathcal{Z} \cap
G$ is, in general, a two dimensional domain.
\end{enumerate}}
\end{rem}

As an example, consider the mapping
\begin{equation}\label{FabEq66}
\psi(w):=\left[\left(w^{-1}-\ou\right)^{1/2}+\left(w^{-1}-\cj{\ou}
\right)^{1/2}+i\,\ou^{1/2}-i\,\cj{\ou}^{1/2}\right]^{-1},\quad
|w|\geq 1,
\end{equation}
where $\ou=e^{ i\Theta_1}$, $ \pi/2 \leq\Theta_1<
\pi$,  is given and the branch of the root chosen
is analytic on $\mathbb{C}\setminus (-\infty,0]$
and positive on $(0,+\infty)$. Then,
$\psi(\infty)=\infty$,
$\psi(w)=\cj{\psi(\cj{w})}$, and $\psi$ maps
$\Delta_1$ conformally onto the exterior of a
piecewise analytic Jordan curve $L$ symmetric
about the real axis, with corners at
$z_1=\psi(\ou)$, $z_2=\psi(\cj{\ou})$. Here,
\[
\ou_1=\ou,\quad\ou_2=\cj{\ou}=e^{
i(2\pi-\Theta_1)},\quad
\theta_2=1-\Theta_1/\pi,\quad
\lambda_1=\lambda_2=1/2,
\]
\[
\hat{A}_1=-i\,\cj{\ou}\,[\psi(\ou)]^2,\quad
\hat{A}_2=i\,[\psi(\cj{\ou})]^2.
\]

\begin{figure}[!h]
\centering
\includegraphics[scale=.55]{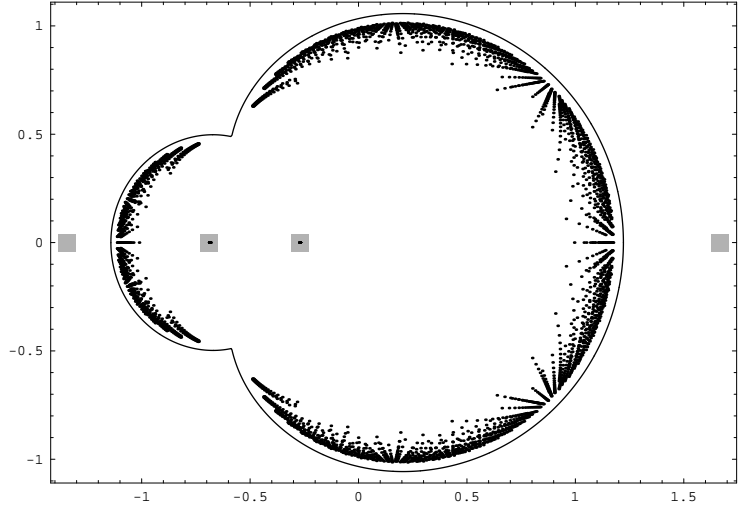}
\caption{Zeros of $F_n(z)$, $n=\cj{20,90}$, for a
$\psi$ as in (\ref{FabEq66}) with $\ou=\exp(3\pi
i/4)$.}\label{FabFig2}
\end{figure}

Therefore, when $\theta_2=p_2/q_2$ is rational,
$1\leq p_2<q_2$ relatively prime integers, a
point $t$ interior to $L$ is an accumulation
point of the zeros of the $F_n$'s if and only if
$t$ satisfies one of the equations
\begin{equation}\label{FabEq65}
(\hat{A}_1+\hat{A}_2e^{2\pi i
s/q_2})t=\hat{A}_1z_2+\hat{A}_2z_1e^{2\pi i
s/q_2},\quad s=0,1,\ldots,q_2-1.
\end{equation}
Figure \ref{FabFig2} below corresponds to the
case $\Theta_1=3\pi/4$ ($\theta_2=1/4$), where we
have plotted the zeros of $F_n(z)$,
$n=\cj{20,90}$. The solutions of the equations in
(\ref{FabEq65}) are the centers of the grayish
squares.

If $\theta_2$ is irrational, then a point $t$
interior to $L$ is an accumulation point of the
zeros of the $F_n$'s if and only if $t$ is real.
In Figure \ref{FabFig1} below, we have plotted
the zeros of $F_n(z)$ for $n=\cj{20,90}$
corresponding to the case
$\Theta_1=\sqrt{2}\pi/2$.
\begin{figure}[!h]
\centering
\includegraphics[scale=.48]{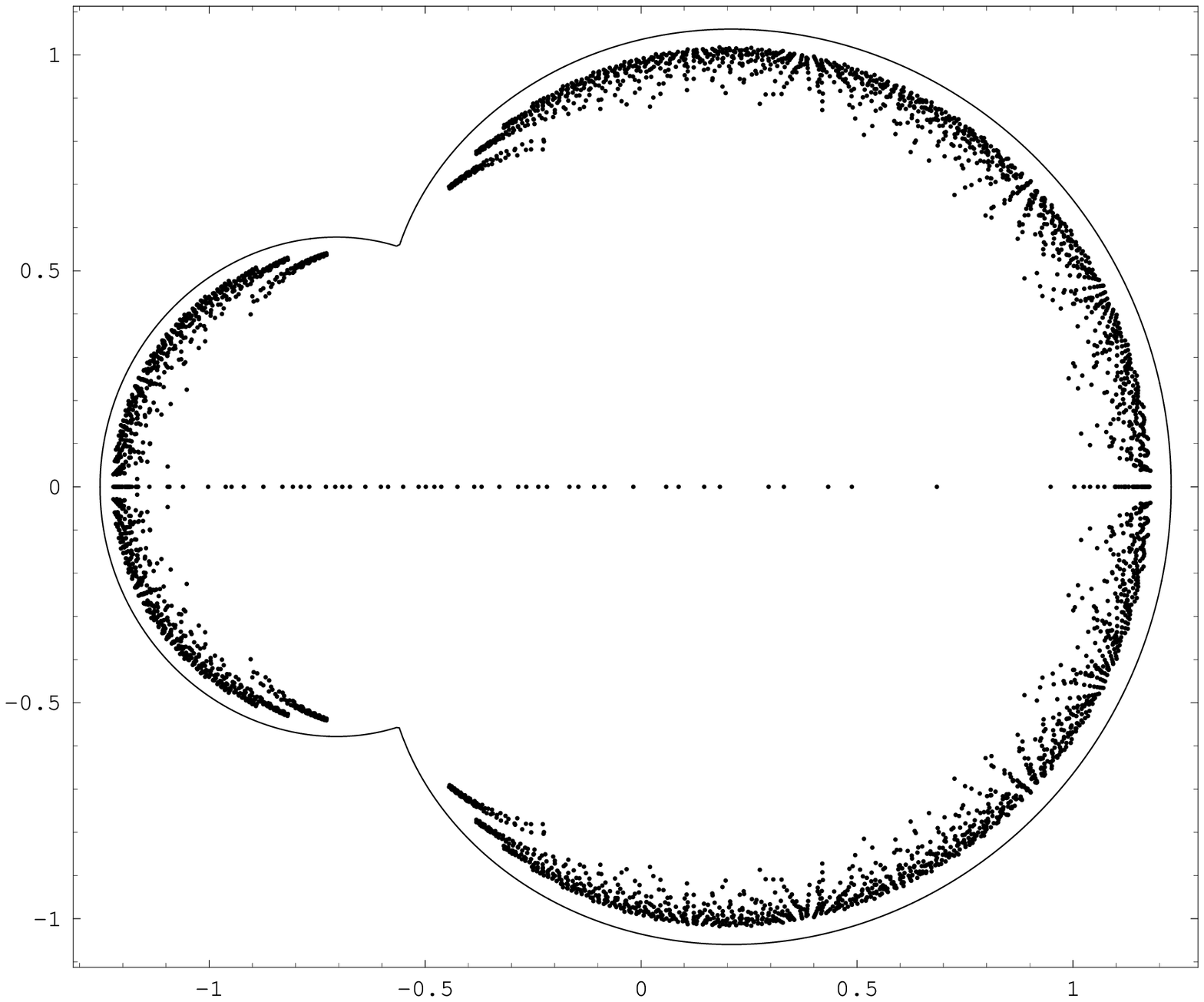}
\caption{Zeros of $F_n(z)$, $n=\cj{20,90}$, for a
$\psi$ as in (\ref{FabEq66}) with
$\ou=\exp(\sqrt{2}\pi i /2)$.}\label{FabFig1}
\end{figure}

\begin{figure}[!h]
\begin{center}
\includegraphics[scale=.5]{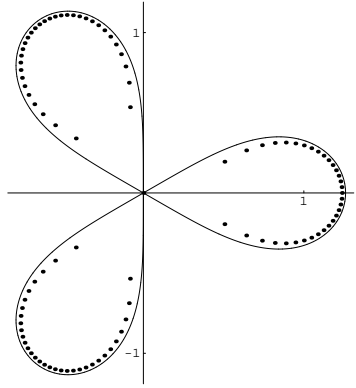}
\end{center}
\caption{Zeros of $F_{100}(z)$ for the domain
$\Omega=\{z:|z^3-1|>1\}$.}\label{FabFig3}
\end{figure}

We finish this section presenting an example in
which condition $A.3$ is not satisfied. Let
$s\geq 2$ be a given integer. The function
$w\mapsto w^s+1$ maps each of the $s$ sectors
\[
2\pi(k-1)/s< \arg(w)<2\pi k/s,\quad
k=1,2,\ldots,s,
\]
conformally onto the complex plane cut along the
ray $[0,+\infty)$, and by agreeing in that
\[
2\pi(k-1)\leq \arg(w^s+1)<2\pi k\quad
\mathrm{whenever}\quad 2\pi(k-1)/s\leq
\arg(w)<2\pi k/s,
\]
we see that
\[
\psi(w)=(w^s+1)^{1/s}
\]
maps the exterior of the unit circle
conformally onto the exterior of the lemniscate
of $s$ petals $L=\{z:|z^s-1|=1\}$ (see Figure
\ref{FabFig3} for $s=3$).

Here, $\Omega=\{z:|z^s-1|>1\}$,
$G=\{z:|z^s-1|<1\}$, and the inverse of $\psi$ is
$\phi(z)=(z^s-1)^{1/s}$. Moreover, it is easily
seen that $\psi$ satisfies conditions $A.1$ and
$A.2$ with
\[
\ou_k=e^{i(2k-1)\pi /s},\quad z_k=0\,, \quad
\lambda_k=1/s,\quad k=1,2,\ldots,s.
\]

The Faber polynomial $F_n(z)$ is the polynomial
part of the Laurent expansion at $\infty$ of
$(z^s-1)^{n/s}$. Hence, for any two integers
$m\geq 0$ and $l\in \{0,1,\ldots,s-1\}$,
\begin{equation}\label{FabEq63}
F_{sm+l}(z)=\sum_{j=0}^m\left(-1\right)^{j}
{m+l/s \choose j} z^{s(m-j)+l},
\end{equation}
where ${ a\choose b}$ stands for the generalized
binomial coefficient
$\Gamma(a+1)/\left[\Gamma(b+1)\Gamma(a-b+1)\right]$.
In particular, $ F_{sm}(z)=\left(z^s-1\right)^m$.

The important feature to note of this example is
that the function $H_n(z)$ defined in
(\ref{FabEq85}) is identically zero for every
$n\not=s-1\mod s$ (recall Remark \ref{FabRem2}).
This example has been previously studied by
Ullman \cite{Ullman} for $s=2$, and by He
\cite{He2} for $s\geq 2$. Observe from
(\ref{FabEq63}) that $F_{sm+l}(z)$ has a zero of
multiplicity $l$ at the origin. Ullman and He
showed that all other zeros  lie strictly in $G$
(see Figure \ref{FabFig3} below).

Theorem \ref{FabThm3} below shows that for every
$l\in \{1,2,\ldots,s-1\}$, we can properly
normalize the subsequence
$\left\{F_{sm+l}\right\}_{m\geq 0}$ so as to make
it converge locally uniformly on $G$ to a
function that never vanishes on $G$. Hence, for
every compact set $E\subset G$, there exists
$N_E$ such that if $n\not=0\mod s$, and $n>N_E$,
then $F_n$ has no zeros on $E$. As a consequence,
$\nu_{n}\wc \mu_L$ as $n\to\infty$, $n\not=0\mod
s$, where
\begin{equation*}
\mu_L=\frac{|z|^{s-1}|dz|}{2\pi},\quad z\in L,
\end{equation*}
is the equilibrium measure of $L$. Observe how
the distribution function of $\mu_L$ is in total
agreement with the density pattern followed by
the zeros in Figure \ref{FabFig3}.

However, the zeros of $F_{sm}$ are fixed, namely
$e^{2\pi ik/s}$, $1\leq k\leq s$, each of
multiplicity $m$ and contained in $G$, and so
\[
\nu_{sm}\wc \frac{1}{s}\sum_{k=1}^s\delta_{e^{2\pi
ik/s}}\quad \mathrm{as}\ m\to\infty.
\]
Thus, Corollary \ref{FabPro1} does not
necessarily hold in the absence of condition A.3.

\begin{thm} \label{FabThm3} For every $l\in
\{1,2,\ldots,s-1\}$,
\begin{equation*}
(-1)^{m+1}{sm+l \choose l/s-1 }^{-1}F_{sm+l}(z)=
\frac{1}{s^{l/s} z^{s-l}}\left[1+r_m(z)\right]
\end{equation*}
where
\[
r_m(z)=\frac{1}{m}\left[\frac{(s-1)(s-l)(2s-l)}{2s^3}-\frac{1}{sz^s}\right]+\mathcal{O}\left(m^{-(1+s)/s}\right)
\]
locally uniformly in $|z^s-1|<1$ as $m\to\infty$.
\end{thm}

More important than Theorem \ref{FabThm3} is its proof,
which illustrates an approach to obtaining asymptotics for
$F_n$ in cases where A.3 is not satisfied.

\section{Lehman expansion of $\psi$ near $\ou_k$}\label{FabSec1}  Let $\zeta$ be a small open
circular arc of $\mathbb{T}_1$ centered at
$\ou_k$ such that
$\cj{\zeta}\cap\{\ou_1,\ldots,\ou_s\}=\{\ou_k\}$.
The set $\zeta\setminus\{\ou_k\}$ consists of two
circular arcs, say $\zeta^+$, $\zeta^-$, and by
our assumption A.1 on $L$, there exist simple
analytic arcs
$\mathcal{L}^+\supset\psi\left(\cj{\zeta^+}\right)$
and
$\mathcal{L}^-\supset\psi\left(\cj{\zeta^-}\right)$
of which $z_k$ is an interior point. Hence the
map $\psi$, originally defined on $\Delta_1$, can
be continued by the Schwarz reflection principle
for analytic arcs \cite{Davis} across both
$\zeta^+$ and $\zeta^-$. Since the images of
$\mathcal{L}^+$ and $\mathcal{L}^-$ in such
reflections are again simple analytic arcs
containing $z_k$ as an interior point, by
applying subsequent reflections we can continue
$\psi$ near $\ou_k$ onto the entire logarithmic
Riemann surface $\mathcal{S}_{\ou_k}$ with branch
point at $\ou_k$.

Let the functions $(w-\ou_k)^{l+j\lambda_k}$,
$l\geq 0$, $j\geq 1$, and $\log(w-\ou_k)$ be
defined in $\mathcal{S}_{\ou_k}$. In what follows
we abbreviate by putting $y=w-\ou_k$. Lehman
\cite[Thm. 1]{Lehman} proved that $\psi$ has the
following asymptotic expansion: if $\lambda_k$ is
\emph{irrational}, then
\begin{equation}\label{FabEq1}
\psi(w)=\psi(\ou_k)+\sum_{l= 0}^\infty \sum_{j=1}^\infty
c^{k}_{lj0}y^{l+j\lambda_k} ,\quad c^{k}_{010}\not=0\,;
\end{equation}
if $\lambda_k=p/q$ is a \emph{fraction} reduced to lowest terms,
then
\begin{eqnarray}\label{FabEq2}
\psi(w)&=&\psi(\ou_k)+\sum_{l=0}^\infty \sum_{j=1}^q
\sum_{m=0}^{\lfloor l/p\rfloor}c^{k}_{ljm}y^{l+j\lambda_k}(\log
y)^m,\quad c^{k}_{010}\not=0.
\end{eqnarray}
The terms in the above series are assumed to be arranged in an order
such that a term of the form $y^{l+j\lambda_k}(\log y)^m$ precedes
one of the form $y^{l'+j'\lambda_k}(\log y)^{m'}$ if either
$l+j\lambda_k<l'+j'\lambda_k$ or $l+j\lambda_k=l'+j'\lambda_k$ and
$m>m'$.

The precise meaning of these expansions is the
following: if according to the order explained
above, (\ref{FabEq1}) and (\ref{FabEq2}) are
written in the form
\[
\psi(w)=\psi(\ou_k)+\sum_{n=1}^\infty \chi_n(y),
\]
then for all $N\geq 1$,
\[
\psi(w)-\psi(\ou_k)-\sum_{n=1}^N \chi_n(y)=o\left(\chi_N(y)\right)
\]
as $w\to \ou_k$ from any finite sector
$\vartheta_1\leq\arg(w-\ou_k)\leq \vartheta_2$ of
$\mathcal{S}_{\ou_k}$.

We write in (\ref{FabEq1}) $c^k_{lj0}$ instead of
simply $c^k_{lj}$ when $\lambda_k$ is irrational,
because this will allow us to express many of the
relations that follow in one single statement
without having to distinguish between $\lambda_k$
being irrational or rational.

The coefficients $c^k_{ljm}$ in (\ref{FabEq1})
and (\ref{FabEq2}) depend on the values assigned
to the functions $(w-\ou_k)^{l+j\lambda_k}$,
$\log(w-\ou_k)$ at a specified point of
$\mathcal{S}_{\ou_k}$. We shall assume that the
values of $\psi$ in $\Delta_1$  define $\psi$ in
the sector
$\Theta_k-\pi<\arg(w-\ou_k)<\Theta_k+\pi$ of
$\mathcal{S}_{\ou_k}$, and that for every $w$ in
this sector,
\[
(y)^{l+j\lambda_k}=|y|^{l+j\lambda_k}e^{i(l+j\lambda_k)\arg(y)},\quad
\log y=\log|y|+i\arg(y), \quad y=w-\ou_k.
\]

A more detailed description of these expansions is split in
two cases:

\noindent \textbf{Case
{\mathversion{bold}$0<\lambda_k<2$,
$\lambda_k\not\in\{1,2\}$}:} As in Section
\ref{FabSubSec1}, we put
$A_k:=c^{k}_{010}\not=0$, and it follows from
(\ref{FabEq1}) and (\ref{FabEq2}) that for
$\upsilon>0$ sufficiently small, say
\[
\upsilon <\left\{\begin{array}{cc}
                   \min\{\lambda_k,1-\lambda_k\}, &\quad\mathrm{if}\ 0<\lambda_k<1, \\
                   2-\lambda_k\,, &\quad\mathrm{if}\ 1<\lambda_k<2,
                 \end{array}
\right.
\]
the following relations hold: \emph{if} $0<\lambda_k<1$, then
\begin{equation}\label{FabEq51}
\psi(w)=\psi(\ou_k)+A_k y^{\lambda_k}+c^{k}_{020}y^{2\lambda_k}+
o\left(y^{2\lambda_k+\upsilon}\right);
\end{equation}
\emph{if} $\lambda_k=1/2$,
\begin{equation}\label{FabEq6}
\psi(w)=\psi(\ou_k)+A_k
y^{\lambda_k}+c^{k}_{020}y^{2\lambda_k}+c^{k}_{111}y^{1+\lambda_k}(\log
y)+
c^{k}_{110}y^{1+\lambda_k}+o\left(y^{3\lambda_k+\upsilon}\right);
\end{equation}
\emph{if} $1<\lambda_k<2$,
\begin{equation}\label{FabEq52}
\psi(w)=\psi(\ou_k)+A_k
y^{\lambda_k}+c^{k}_{110}y^{1+\lambda_k}+c^{k}_{020}y^{2\lambda_k}+
o\left(y^{2\lambda_k+\upsilon}\right)
\end{equation}
(notice that if $1<\lambda_k=p/q<2$, then $p\geq
3$, $q\geq 2$, and no $\log$-terms correspond to
$l=0,1,2$).

\noindent \textbf{Case{\mathversion{bold}
$\lambda_k\in\{1,2\}$}:} Here $p=\lambda_k,\
q=1$. \emph{If $\ou_k$ is relevant}, then there
is a smallest integer $r_k\geq \lambda_k$ for
which a $\log$-term of the form
$y^{r_k+\lambda_k}(\log y)^{m_k}$, $1\leq m_k\leq
\lfloor r_k/\lambda_k\rfloor$, occurs in the
expansion of $\psi$ about $\ou_k$, so that in
case \emph{$m_k\geq 2$},
\begin{eqnarray}\label{FabEq10}
\psi(w)&=&\psi(\ou_k)+\sum_{l=0}^{r_k-1}
c^{k}_{l10}y^{l+\lambda_k}+A_{k}y^{\Lambda_k}(\log
y)^{m_k}+B_{k}y^{\Lambda_k}(\log y)^{m_k-1}\\
&&+C_k y^{\Lambda_k}(\log y)^{m_k-2}+\left\{\begin{array}{ll}
\underset{\,}{\mathcal{O}\left(y^{\Lambda_k}(\log y)^{m_k-3}\right)}, &\quad \mathrm{if}\ m_k\geq 3,\\
\mathcal{O}\left( y^{\Lambda_k+1/2} \right),&
\quad \mathrm{if}\ m_k= 2,
\end{array}
\right.\nonumber
\end{eqnarray}
while \emph{if $m_k=1$}, then
\begin{eqnarray}\label{FabEq17}
\psi(w)&=&\psi(\ou_k)+\sum_{l=0}^{r_k-1}
c^{k}_{l10}y^{l+\lambda_k}+A_{k}y^{\Lambda_k}(\log
y)^{m_k}+B_{k}y^{\Lambda_k}(\log
y)^{m_k-1}\nonumber\\
&&+\widetilde{C}_k y^{\Lambda_k+1}(\log y)^{\lfloor(r_k+1)/\lambda_k
\rfloor}+
\widetilde{D}_k y^{\Lambda_k+1}(\log y)^{\lfloor(r_k+1)/\lambda_k\rfloor-1}\nonumber\\
&&+               \left\{\begin{array}{ll}
                          \mathcal{O}\left( y^{\Lambda_k+1}(\log
y)^{\lfloor(r_k+1)/\lambda_k \rfloor-2}\right), &
\quad\mathrm{if}\
\lfloor(r_k+1)/\lambda_k\rfloor\geq 2,\\
\mathcal{O}\left( y^{\Lambda_k+3/2}\right), &
\quad\mathrm{if}\
\lfloor(r_k+1)/\lambda_k\rfloor=1,
                        \end{array}\right.
\end{eqnarray}
where
\begin{equation}
\Lambda_k:=r_k+\lambda_k, \ \, A_k:=c^{k}_{r_k1m_k}\not=0,\ \,
B_k:=c^{k}_{r_k1(m_k-1)},\ \, C_k:=c^{k}_{r_k 1(m_k-2)},
\end{equation}
\begin{equation}
\widetilde{C}_k:=c^{k}_{(r_k+1)1\lfloor(r_k+1)/\lambda_k \rfloor},
\quad \widetilde{D}_k:=c^{k}_{(r_k+1)1\left(\lfloor(r_k+1)/\lambda_k
\rfloor-1\right)}\,.
\end{equation}
Thus, setting
\begin{equation}\label{FabEq14}
Q_k(w):=z_k+\sum_{l=0}^{r_k-1}
c^{k}_{l10}y^{l+\lambda_k}=z_k+\mathcal{O}\left(y^{\lambda_k}\right),
\end{equation}
we have that \emph{if $m_k\geq 2$}, then
\begin{eqnarray}\label{FabEq13}
\psi'(w)&=&Q_k'(w)+A_{k}\Lambda_k y^{\Lambda_k-1}(\log
y)^{m_k}+\left(A_k m_k+B_k \Lambda_k\right)y^{\Lambda_k-1}(\log y)^{m_k-1}\nonumber\\
&&+\left[B_k(m_k-1)+C_k\Lambda_k\right]
y^{\Lambda_k-1}(\log
y)^{m_k-2}+\left\{\begin{array}{ll}
                \mathcal{O}\left(y^{\Lambda_k-1}(\log y)^{m_k-3}\right), &\quad \mathrm{if}\ m_k\geq 3,\\
               \mathcal{O}\left( y^{\Lambda_k-1/2} \right),& \quad \mathrm{if}\ m_k=
               2,
              \end{array}
\right.\nonumber\\
\end{eqnarray}
while \emph{if $m_k=1$}, then
\begin{eqnarray}\label{FabEq12}
\psi'(w)&=&Q_k'(w)+A_{k}\Lambda_k y^{\Lambda_k-1}(\log
y)^{m_k}+\left(A_k m_k+B_k \Lambda_k\right)y^{\Lambda_k-1}(\log y)^{m_k-1}\nonumber\\
&&+\widetilde{C}_k\left(\Lambda_k+1\right) y^{\Lambda_k}(\log
y)^{\lfloor(r_k+1)/\lambda_k
\rfloor}\nonumber\\
&&+\left[\widetilde{C}_k\lfloor(r_k+1)/\lambda_k \rfloor+
\widetilde{D}_k \left(\Lambda_k+1\right)\right]y^{\Lambda_k}(\log y)^{\lfloor(r_k+1)/\lambda_k\rfloor-1}\nonumber\\
&&+               \left\{\begin{array}{ll}
                          \mathcal{O}\left( y^{\Lambda_k+1/2}\right) & \quad\mathrm{if}\
\lfloor(r_k+1)/\lambda_k\rfloor=1,\\
\mathcal{O}\left( y^{\Lambda_k}(\log
y)^{\lfloor(r_k+1)/\lambda_k \rfloor-2}\right) &
\quad\mathrm{if}\
\lfloor(r_k+1)/\lambda_k\rfloor\geq 2\,.
                        \end{array}\right.
\end{eqnarray}

\emph{If $\ou_k$ is not relevant}, then for every
$N\in\mathbb{N}$,
\begin{equation}
\psi(\ou_k)=\psi(\ou_k)+\sum_{l=0}^{N}
c^{k}_{l10}y^{l+\lambda_k}+o\left(y^{N+\lambda_k}\right)=Q_{k}(z)+o\left(y^{N+\lambda_k}\right),
\end{equation}
where
\begin{equation}\label{FabEq53}  Q_{k}(z):=z_k+\sum_{l=0}^{N}
c^{k}_{l10}y^{l+\lambda_k}.
\end{equation}

The polynomial  $Q_{k}(z)$ defined by (\ref{FabEq53})
depends on the value $N$, so that in what follows we will
think of $N$ as an arbitrarily large natural number that
has been fixed.

\section{Proofs of the asymptotic results}\label{FabSec2}

Recall we are using the notation
\[
\mathbb{T}_r:=\{w:|w|=r\},\quad
\mathbb{D}_r:=\{w:|w|<r\},\quad
\Delta_r:=\{w:r<|w|\leq \infty\}.
\]
Also, for $a,\,b\in \mathbb{C}$, we denote by
$[a,b]$ the oriented closed segment that starts
at $a$ and ends at $b$. A similar meaning is
attached to $(a,b)$, $(a,b]$ and $[a,b)$. For
every $0<\sigma<1$, we define
\[
\sigma_k:=\sigma \ou_k\,,\qquad 1\leq k\leq s,
\]
and the contour
\[
\Gamma_\sigma:=\mathbb{T}_\sigma\cup\left(\cup_{k=1}^s[\sigma_k,\ou_k]\right)\,.
\]
The exterior of the contour $\Gamma_\sigma$,
denoted by $\mathrm{ext}(\Gamma_\sigma)$, is
understood to be the unbounded component of
$\cj{\mathbb{C}}\setminus\Gamma_\sigma$, that is,
\[
\mathrm{ext}(\Gamma_\sigma)=\Delta_\sigma\setminus\left(\cup_{k=1}^s[\sigma_k,\ou_k]\right)\,.
\]

Disregarding technical difficulties, the idea
behind the proofs of the asymptotic results is
simple, and in rough terms can be described as
follows. By the piecewise analyticity of
$\partial\Omega$, if $1-\sigma$ is small enough,
the function $w^n\psi'(w)/(\psi(w)-z)$ (in the
variable $w$, for fixed $z$)
 has a meromorphic
extension to $\mathrm{ext}(\Gamma_\sigma)$ with
at most finitely many poles in there and
continuous boundary values on $\Gamma_\sigma$.
Then, using the integral representation
(\ref{FabEq27}), we can express $F_n(z)$ as the
sum of the residues of that function in
$\ext(\Gamma_\sigma)$, plus its integral (with
respect to $w$) over $\Gamma_\sigma$, the later
being split as an integral over
$\mathbb{T}_\sigma$ (which is
$\mathcal{O}(\sigma^n)$ as $n\to\infty$, and
therefore negligible) plus the integral of
$w^n\psi'(w)/(\psi(w)-z)$ over each of the
``two-sided'' segments $[\sigma_k,\ou_k]$, $1\leq
k\leq s$. The asymptotic behavior as $n\to\infty$
of these last integrals (as functions of $z$) can
then be obtained from the Lehman expansions of
$\psi$ about the $\ou_k$'s.

The first step in doing all this rigorously is to prove
that a contour $\Gamma_\sigma$ satisfying the necessary
conditions exists. That is the content of Lemmas
\ref{FabLem2} and \ref{FabLem1} below.

For given $\delta>0$ and $t\in\mathbb{C}$, we put
\[
D_\delta(t):=\{w:|w-t|<\delta\},\quad
D_\delta^*(t):=\{w:0<|w-t|<\delta\}.
\]
Then, for $\delta>0$ sufficiently small, the set
\[
\zeta_{k,\delta}:=D_\delta(\ou_k)\cap\mathbb{T}_1
\]
is a circular arc, and
$\zeta_{k,\delta}\setminus\{\ou_k\}$ is the union
of two disjoint open circular arcs that we denote
by $\zeta_{k,\delta}^+$, $\zeta_{k,\delta}^-$,
say $\zeta_{k,\delta}^-$ immediately followed by
$\zeta_{k,\delta}^+$ when $\mathbb{T}_1$ is
traveled in counterclockwise direction.

Let
\[
D^+_\delta(\ou_k):=D^*_\delta(\ou_k)\setminus
\zeta_{k,\delta}^+,\quad
D^-_\delta(\ou_k):=D^*_\delta(\ou_k)\setminus
\zeta_{k,\delta}^-\,.
\]
If $\delta$ is sufficiently small, the disks
$D_\delta(\ou_1),D_\delta(\ou_2),\dots,D_\delta(\ou_s)$
are pairwise disjoint, and for every
$k\in\{1,2,\ldots,s\}$, the mapping $\psi$ has
analytic continuations $\psi_+$, $\psi_-$ from
the exterior $\Delta_1$ of the unit circle to
$D^+_\delta(\ou_k)$, $D^-_\delta(\ou_k)$,
respectively.

Recall that $L_1$ and $L_2$ are defined by
(\ref{FabEq69})-(\ref{FabEq91}), and that for
every $k$ with $\lambda_k\in\{1,2\}$, $Q_k$ is
defined by (\ref{FabEq14}) and (\ref{FabEq53}).

\begin{lem}\label{FabLem2} Let $\epsilon>0$ be given.  For every $\delta>0$ sufficiently
small the following statements hold true:
\begin{enumerate}
\item[\emph{(a)}] For all $k\in\{1,2,\ldots,s\}$,
\begin{equation}\label{FabEq70}
\cj{\psi_\pm
\left(D^\pm_{\delta}(\ou_k)\right)}\subset
D_{\epsilon}(z_k)\,
\end{equation}
and
\begin{equation}\label{FabEq93}
\psi_\pm \left(w\right)\not=z_k\,,\quad \psi'_\pm
(w)\not=0\quad \forall\, w\in
D^\pm_{\delta}(\ou_k).
\end{equation}
Also, for every $k$ with $\lambda_k\in
\{1,2\}$, we have that
\begin{equation}\label{FabEq94}
Q_k(w)\not=z_k\quad \forall\, w\in
D^*_{\delta}(\ou_k).
\end{equation}

\item[\emph{(b)}] If $\tau_k:=(1-\delta)\ou_k$, there is a
constant $C$ such that for every $k$ with $z_k\in L_1$
(resp. $z_k\in L_2$),
\begin{equation}\label{FabEq31}
\left|\frac{\psi_\pm(t)-z_k}{\psi_\pm(t)-z}\right|<
C,\quad
\left|\frac{Q_k(t)-z_k}{Q_k(t)-z}\right|<C,
\end{equation}
for all $t\in (\tau_k,\ou_k)$ and $z\in
\cj{\Omega}$ (resp. $z\in L_2$).

\end{enumerate}
\end{lem}
\begin{demo}{\emph{\textbf{Proof of Lemma \ref{FabLem2}}}} Since for every $k\in\{1,2,\ldots,s\}$, $\psi$
has an analytic continuation to the entire
logarithmic Riemann surface with branch point at
$\ou_k$, and it is such that
\begin{equation}\label{FabEq61}
\psi(w)-z_k=c^{k}_{010}(w-\ou_k)^{\lambda_k}(1+o(1))\,,\quad
c^{k}_{010}\not=0,
\end{equation}
\begin{equation}\label{FabEq95}
\psi'(w)=c^{k}_{010}\lambda_k(w-\ou_k)^{\lambda_k-1}(1+o(1))\,,
\end{equation}
as $w\to \ou_k$ from any finite sector of the surface, and
since by the very definition of $Q_k$ in (\ref{FabEq14})
and (\ref{FabEq53}), for every $k$ with $\lambda_k\in
\{1,2\}$
\begin{equation}\label{FabEq62}
Q_k(w)-z_k=c^{k}_{010}(w-\ou_k)^{\lambda_k}(1+o(1)),
\end{equation}
as $w\to \ou_k$, the conditions of part (a) of
Lemma \ref{FabLem2} will be trivially satisfied
provided that $\delta$ is small enough.

Let us now prove part (b) of the lemma. The
analysis is split in two cases.

\noindent\emph{\underline{Case 1}: $k$ is such
that $z_k\in L_1\cup L_2$ and $\lambda_k=1$.}
Under this assumption, there is a small open
circular arc $\zeta_k\subset\mathbb{T}_1$, with
center point $\ou_k$, such that $\psi$ is
one-to-one on $\zeta_k$ and therefore
$\psi(\zeta_k)$ is a simple smooth arc containing
$z_k$ as an interior point. We first verify the
following

 \emph{Claim: There is a small open disk
$O_k$ centered at $z_k$ such that if $z_k\in
L_2$, then $O_k\cap L_2\subset \psi(\zeta_k)$,
while if $z_k\in L_1$, then $O_k\cap L\subset
\psi(\zeta_k)$ and $O_k\setminus \psi(\zeta_k)$
is the disjoint union of two nonempty connected
open sets, one contained in $\Omega$, the other
in $G$.}

Let us prove the claim. Suppose $z_k\in L_2$ and
let $t_k\not=\ou_k$ be the other point of
$\mathbb{T}_1$ such that $\psi(t_k)=z_k$. Then,
there is also a small open circular arc
$\rho_k\subset\mathbb{T}_1$, with center point
$t_k$, such that $\psi$ is one-to-one on
$\rho_k$. Since $\eta(z_k)=2$, the closed set
$\psi\left(\mathbb{T}_1\setminus(\zeta_k\cup
\rho_k)\right)$ cannot contain $z_k$, so that
there is a small open disk $O_k$ centered at
$z_k$ such that $\{t\in \mathbb{T}_1:\psi(t)\in
O_k\}\subset (\zeta_k\cup \rho_k)$. But since
$\psi$ is one-to-one on $\rho_k$, we must have
that $O_k\cap L_2\subset \psi(\zeta_k)$.
Similarly, suppose $z_k\in L_1$ and let
$\delta'>0$ be so small that the connected open
set $\psi(D_{\delta'}(\ou_k)\cap\Delta_1)$ lies
strictly on one side of the arc $\psi(\zeta_k)$.
Since $\eta(z_k)=1$, the closed set
$\psi\left(\cj{\Delta}_1\setminus
D_{\delta'}(\ou_k)\right)$ cannot contain $z_k$,
and therefore, there is a sufficiently small open
disk $O_k$ centered at $z_k$ such that $O_k\cap
\psi\left(\cj{\Delta}_1\setminus
D_{\delta'}(\ou_k)\right)=\emptyset$. In
consequence, $O_k$ is divided by the arc
$\psi(\zeta_k)$ into two connected open sets, one
contained in
$\psi(D_{\delta'}(\ou_k)\cap\Delta_1)\subset
\Omega$, the other contained in $G$. The claim is
proven.

Now, for every $k$ with $z_k\in L_1\cup L_2$,
choose $O_k$ as in the claim, and assume
$\delta>0$ is so small that, besides satisfying
part (a) of the lemma, it also satisfies that
$\cj{\psi_\pm
\left(D^\pm_{\delta}(\ou_k)\right)}\subset O_k$.
By our assumption A.1 on $L$, there are two
simple analytic arcs $ \mathcal{L}_k^+$,
$\mathcal{L}_k^-$, each containing $z_k$ as an
interior point, and such that
$\psi\left(\cj{\zeta_{k,\delta}^+}\right)\subset
\mathcal{L}_k^+$,
$\psi\left(\cj{\zeta_{k,\delta}^-}\right)\subset
\mathcal{L}_k^-$. Notice that $\mathcal{L}^\pm_k$
and $\psi(\zeta_k)$ share the same tangent line
at $z_k$. In consequence, if
$\tau_k:=(1-\delta)\ou_k$, the arc
$\psi\left(\left(\ou_k,1/\cj{\tau}_k\right)\right)$
lies entirely in $\Omega$, and by
(\ref{FabEq61}), it is perpendicular to
$\psi(\zeta_k)$.

By the Schwarz reflection principle for analytic
arcs \cite{Davis}, if $\tau_k$ is close enough to
$\ou_k$, $\psi_\pm((\tau_k,\ou_k))$ is the
reflection of
$\psi\left(\left(\ou_k,1/\cj{\tau}_k\right)\right)$
across $\mathcal{L}^\mp_k$, and therefore,
\emph{for all $\delta$ sufficiently small, the
arc $\psi_\pm((\tau_k,\ou_k))$ is perpendicular
to $\psi(\zeta_k)$ and}
\begin{equation}\label{FabeEq63}
\psi_\pm\left((\tau_k,\ou_k)\right)\subset
G,\quad \psi\left(\mathbb{D}_1\cap
D_{\delta}(\ou_k)\setminus
(\tau_k,\ou_k]\right)\subset G\quad if\ z_k\in
L_1,
\end{equation}
\begin{equation}\label{FabEq90}
\psi_\pm\left((\tau_k,\ou_k)\right)\cap
L_2=\emptyset ,\quad\psi\left(\mathbb{D}_1\cap
D_{\delta}(\ou_k)\setminus
(\tau_k,\ou_k]\right)\cap L_2=\emptyset \quad if\
z_k\in L_2,\ \lambda_k=1,
\end{equation}
whence it follows at once that the first
inequality of (\ref{FabEq31}) holds true.

Similar considerations apply to $Q_k$. Because of
(\ref{FabEq62}), if $z_k\in L_1\cup L_2$ and
$\lambda_k=1$, then $Q_k$ maps a small circular
arc of $\mathbb{T}_1$ centered at $\ou_k$ onto an
analytic arc tangent to $ \psi(\zeta_k)$ at
$z_k$, and therefore, for all $\delta$
sufficiently small, $Q_k((\tau_k,\ou_k))$ is
perpendicular to $\psi(\zeta_k)$ and is contained
entirely in $G$, whence the second inequality of
(\ref{FabEq31}) follows.

\emph{\underline{Case 2}: $k$ is such that
$\lambda_k=2$.} Let $\zeta_k\subset\mathbb{T}_1$
be a small open circular arc with center point at
$\ou_k$, so that $\zeta_k\setminus \{\ou_k\}$
splits into two disjoint circular arcs
$\zeta^+_k$ and $\zeta^-_k$. If $\zeta_k$ is
small enough, $\psi\left(\cj{\zeta^+_k}\right)$
and $\psi\left(\cj{\zeta^-_k}\right)$ are simple
analytic arcs forming a cusp pointing toward
$\Omega$. Since $\lambda_k=2$ and all
$\lambda_j$'s are strictly positive (that is,
$\Omega$ has not outward pointing cusps),
$\eta(z_k)=1$ and there is an open disk $O_k$
centered a $z_k$ such that $O_k\cap L\subset
\psi(\zeta_k)$. But if two analytic arcs coincide
at infinitely many points, they must be part of
one and the same arc, so that if $O_k$ is
sufficiently small, either $O_k\cap L_2=\{z_k\}$
or $z_k$ is the endpoint of a cut, that is, one
of the two arcs
$\psi\left(\cj{\zeta^+_k}\right)$,
$\psi\left(\cj{\zeta^-_k}\right)$ is contained in
the other.

Reasoning as we did for the case $\lambda_k=1$
above, we derive from the Schwarz reflection
principle that for all $\delta$ sufficiently
small, the arc $\psi_\pm((\tau_k,\ou_k))$ forms
angle $\pi$ with each of the arcs
$\psi\left(\cj{\zeta^+_k}\right)$,
$\psi\left(\cj{\zeta^-_k}\right)$ and
\begin{equation}\label{FabEq103}
\psi_\pm\left((\tau_k,\ou_k)\right)\cap
L_2=\emptyset ,\quad\psi\left(\mathbb{D}_1\cap
D_{\delta}(\ou_k)\setminus
(\tau_k,\ou_k]\right)\cap L_2=\emptyset,
\end{equation}
whence the first inequality of (\ref{FabEq31})
easily follows.

Similarly, by (\ref{FabEq62}), if $\lambda_k=2$,
then for all $\delta$ sufficiently small,
$Q_k((\tau_k,\ou_k))$ forms angle $\pi$ with each
of the arcs $\psi\left(\cj{\zeta^+_k}\right)$,
$\psi\left(\cj{\zeta^-_k}\right)$, whence the
second inequality of (\ref{FabEq31}) follows.

\end{demo}
\begin{lem}\label{FabLem1}
Let $E$ be a (fixed) closed set $(\infty\not\in
E)$ such that either $E\subset G$, or
$E\subset\Omega\cup L_1$ or $E\subset L_2$, or
$E=\{z_1,\ldots,z_s\}$. Then for all $\sigma<1$
with $1-\sigma$ sufficiently small, we have
\begin{enumerate}
\item[\emph{(a)}]  $\psi$ has an analytic continuation to
$\mathrm{ext}(\Gamma_\sigma)$ with
$\psi'(w)\not=0$ for all $w\in
\mathrm{ext}(\Gamma_\sigma)$, and both $\psi$ and
$\psi'$ have continuous boundary values on
$\Gamma_\sigma\setminus\{\ou_1,\ldots\ou_s\}$
when viewing each $[\sigma_k,\ou_k]$ as having
two sides;

\item[\emph{(b)}] if $E\subset G$, then for every $z\in E$, $\psi'(w)/(\psi(w)-z)$
is analytic on $\ext(\Gamma_\sigma)$ with
continuous boundary values on
$\Gamma_\sigma\setminus\{\ou_1,\ldots,\ou_s\}$;

\item[\emph{(c)}] if either $E\subset\Omega\cup L_1$, or $E\subset
L_2$, or $E=\{z_1,\ldots,z_s\}$, then for every
$z\in E$, $\psi'(w)/(\psi(w)-z)$ is analytic on
\[\ext(\Gamma_\sigma)\setminus\left\{\phi_1(z),\ldots,\phi_{\eta(z)}(z)\right\},\]
with a simple pole at each $\phi_j(z)\not\in
\{\ou_1,\ldots,\ou_s\}$ and continuous boundary
values on
$\Gamma_\sigma\setminus\{\ou_1,\ldots,\ou_s\}$ .
\end{enumerate}
\end{lem}

\begin{demo}{\emph{\textbf{Proof of Lemma \ref{FabLem1}}}}

\emph{Part (a) and (b):} Let $\delta>0$ be such
that for every $k\in\{1,2,\ldots,s\}$, the
analytic continuations $\psi_{\pm}$ of $\psi$ to
$D^{\pm}_{\delta}(\ou_k)$ satisfy that
\begin{equation}\label{FabEq104}
\psi'_{\pm}(w)\not=0\quad \forall\, w\in
D^{\pm}_{\delta}(\ou_k).
\end{equation}
Fix $\delta'$ with $0<\delta'<\delta$, and for
every $\sigma$ with $1-\delta'<\sigma<1$,
consider the open set
\[
A_{\delta',\sigma}:=\left\{w:\sigma<|w|<1/\sigma,\
w\not\in\cup_{k=1}^s D_{\delta'}(\ou_k)\right\},
\]
that consists of $s$ open components
$A^l_{\delta',\sigma}$, $l=1,2,\ldots,s$. Then,
by assumption A.1 on $L$, the univalency of
$\psi$ on $\Delta_1$, and the way analytic
functions are continued across analytic arcs by
means of the Schwarz reflection principle, we
have that if $1-\sigma$ is small enough, then
$\psi$ has an analytic and univalent continuation
to each $A^l_{\delta',\sigma}$. From this and
(\ref{FabEq104}), it follows that statement (a)
holds for all $\Gamma_\sigma$ with $\sigma$
sufficiently close to $1$. Moreover, if $E\subset
G$ and $\sigma$ is so close to $1$ that
$\cj{\psi(\mathrm{ext}(\Gamma_\sigma))}\cap
E=\emptyset$, then (b) obviously holds.

\emph{Part (c):} Suppose first that either
$E\subset\Omega\cup L_1$ or $E\subset L_2$. Let
$\epsilon$ be such that\footnote{$
\mathrm{dist}(A,B)=\inf\{|x-y|:(x,y)\in A\times
B\}$.}
\[
0<\epsilon<\mathrm{dist}\left(E,\left\{z_k:z_k\not\in
E,\ 1\leq k\leq s \right\}\right)
\]
For this $\epsilon$, choose $\delta>0$ for which
(\ref{FabeEq63}), (\ref{FabEq90}) and
(\ref{FabEq103}) hold true.

Let $E_1^{-1}:=\left\{w\in
\mathbb{T}_1\setminus\cup_{k=1}^s
D_{\delta}(\ou_k):\psi(w)\in E\right\}$. Then
$E_1^{-1}$ is a compact set, and again, by the
Schwarz reflection principle, we can find a
finite set of open disks $U_1,U_2,\ldots,U_m$,
each centered at some point of $E_1^{-1}$, such
that $E_1^{-1}\subset \cup_{j=1}^m U_j$ and for
all $1\leq j\leq m$, $\psi$ has an analytic and
univalent continuation to each $U_j$, which
satisfies
\begin{enumerate}
\item[i)] $\psi\left(\mathbb{D}_1\cap U_j\right)\subset
G$, in case $E\subset \Omega\cup L_1$,
\item[ii)]$\psi\left(U_j\setminus\mathbb{T}_1\right)\cap
L_2=\emptyset$, in case $E\subset L_2$.
\end{enumerate}

On the other hand, if $E=\{z_1,\ldots,z_s\}$,
choose $\epsilon>0$ such that
\begin{equation}\label{FabEq106}
D_{\epsilon}(z_k)\cap
D_{\epsilon}(z_j)=\emptyset\quad\mathrm{
whenever}\quad z_k\not=z_j,
\end{equation}
and for this $\epsilon$ choose $\delta>0$ so that
Lemma \ref{FabLem2} holds. In this case,
\[E_1^{-1}:=\left\{w\in
\mathbb{T}_1\setminus\cup_{k=1}^s
D_{\delta}(\ou_k):\psi(w)\in E\right\}\] has,
say, $m\geq 0$ elements and we can find $m$ open
disks $U_1,U_2,\ldots,U_m$, each centered at some
$w\in E^{-1}_1$,  such that for all
$j\in\{1,\ldots,m\}$, $\psi$ has an analytic and
univalent continuation to each $U_j$, which
satisfies that
\begin{enumerate}
\item[iii)] $\psi\left(U_j\setminus \mathbb{T}_1\right)\cap
\{z_1,\ldots,z_s\}=\emptyset$.
\end{enumerate}

Now, in either of the three cases i), ii) and
iii) above, the set
 \[
 Y:=\mathbb{T}_1\setminus\left[\left(\cup_{k=1}^s D_{\delta}(\ou_k)\right)\cup\left(\cup_{j=1}^m U_j\right)
\right]
\]
is compact, $\psi(Y)\cap E=\emptyset$, and
therefore there is a neighborhood $W_Y$ of $Y$
such that the analytic continuation of $\psi$ to
$W_Y$ satisfies
\begin{enumerate}
\item[iv)] $\psi(W_Y)\cap E=\emptyset$.
\end{enumerate}
Then, part (c) holds for every $\sigma$ so close
to 1 that
\[
\left\{w:\sigma\leq|w|<1\right\}\subset
W_Y\cup\left(\cup_{k=1}^sD_{\delta}(\ou_k)\right)\cup\left(\cup_{j=1}^m
U_j\right).
\]
This follows for $E\subset\Omega\cup L_1$ from
(\ref{FabeEq63}), i) and iv); for $E\subset L_2$
from (\ref{FabEq90})-(\ref{FabEq103}), ii) and
iv); and for $E=\{z_1,\ldots,z_s\}$ from Lemma
\ref{FabLem2}(a), (\ref{FabEq106}),  iii) and
iv).
\end{demo}

\begin{lem}\label{FabLem3} Let $\epsilon>0$ and $\delta>0$ be such that $\cj{\psi_\pm
\left(D^\pm_{\delta}(\ou_k)\right)}\subset
D_\epsilon(z_k)$. Then, for every $\sigma$ with
$0<1-\sigma<\delta$, we have
that\begin{enumerate}
\item[\emph{i)}] if $\ou_k$ is relevant (that is, if $1\leq k\leq v$), then
\begin{equation}\label{FabEq54}
\frac{1}{2\pi
i}\int_{\sigma_k}^{\ou_k}\left(\frac{t^n\psi_+'(t)}{\psi_+(t)-z}-\frac{t^n\psi_-'(t)}{\psi_-(t)-z}\right)dt=
\alpha_{\Lambda_k-1,M_k}(n)
\left(\frac{\mathcal{C}_kA_k
e^{i(n+\Lambda_k)\Theta_k}}{z-z_k}
+r_{\sigma_k,n}(z)\right)
\end{equation}
with $r_{\sigma_k,n}(z)$ converging uniformly to
zero on $\{z:|z-z_k|\geq\epsilon\}$ as
$n\to\infty$, and
\[
\mathcal{C}_k:=\left\{\begin{array}{ll}
                        [\Gamma(-\lambda_k)\Gamma(\lambda_k)]^{-1}, &\quad if\ \lambda_k\not\in\{1,2\}, \\
                        (-1)^{\Lambda_k-1}m_k\Lambda_k, & \quad if\
                        \lambda_k\in\{1,2\};
                      \end{array}
\right.
\]

\item[\emph{ii)}] if $\ou_k$ is not relevant, then for every $\tau>0$,
\[
\frac{1}{2\pi
i}\int_{\sigma_k}^{\ou_k}\left(\frac{t^n\psi_+'(t)}{\psi_+(t)-z}-\frac{t^n\psi_-'(t)}
{\psi_-(t)-z}\right)dt
=\mathcal{O}\left(n^{-\tau}\right)
\]
uniformly on $\{z:|z-z_k|\geq\epsilon\}$ as
$n\to\infty$.
\end{enumerate}
\end{lem}

\begin{rem}\label{FabRem2}\emph{A more detailed version of (\ref{FabEq54}) given by equalities (\ref{FabEq23}),
(\ref{FabEq47}), (\ref{FabEq58}), (\ref{FabEq26})
and (\ref{FabEq25}) in the proof of Lemma
\ref{FabLem3} provides asymptotic formulas for
the functions $r_{\sigma_k,n}(z)$ from which the
following table follows.}\newline

\label{FabPag1}\emph{{\footnotesize \begin{tabular}{|c|c|c|}\hline
if $k$ is such that & rate of
decay of  $r_{\sigma_k,n}(z)$ is & rate is exact iff\\
\hline
 $0<\lambda_k<1,\ \lambda_k\not=1/2$& $\mathcal{O}\left(n^{-\lambda_k}\right)$&$2c^{k}_{020}(z-z_k)+(A_k)^2\not=0$\\
\hline $\lambda_k=1/2$ & $\mathcal{O}\left(n^{-1}\log n\right)$
&$c^{k}_{111}\not=0$  \\ \hline
 $1<\lambda_k<2$ &$\mathcal{O}\left(n^{-1}\right)$ &
$c^{k}_{110}\not=0$
\\ \hline
$\lambda_k\in\{1,2\},\ m_k\geq 2$&$\mathcal{O}\left(1/\log n\right)$
&
$A_km_k+\Lambda_kc^{k}_{r_k1(m_k-1)}+i\Theta_km_kA_k\Lambda_k\not=0$\\
\hline $\lambda_k\in\{1,2\},\ m_k=1$& $\mathcal{O}\left(n^{-1}(\log
n)^{\lfloor(r_k+1)/\lambda_k \rfloor-1}\right)$ &$c^{k}_{(r_k+1)
1\lfloor (r_k+1)/\lambda_k\rfloor}\not=0$ \\\hline
\end{tabular}}}\newline

\emph{Given $k\in\{1,2,\ldots,s\}$ satisfying one
of the conditions listed in the first column of
the table, an estimate on the rate of decay of
$r_{\sigma_k,n}(z)$ holding uniformly as
$n\to\infty$ on any closed set
$E\subset\{z:|z-z_k|\geq \epsilon\}$ is given in
the second column. The rate is exact for given
$k$  and $E$ if and only if the condition in the
third column is satisfied by every $z\in E$.}
\end{rem}

\begin{demo}{\emph{\textbf{Proof of Lemma \ref{FabLem3}}}} Part (a): First, notice that for
every integer $N\geq 1$, we have the identity
\begin{equation}\label{FabEq5}
\frac{1}{\psi(w)-z}=\sum_{j=0}^{N-1}\frac{\left[z_k-\psi(w)\right]^{j}}{\left(z_k-z\right)^{j+1}}
+\frac{\left[z_k-\psi(w)\right]^{N}}{\left(z_k-z\right)^{N}\left[\psi(w)-z\right]}.
\end{equation}

Suppose first that $\lambda_k\not\not\in
\{1,2\}$. Then, combining identity (\ref{FabEq5})
corresponding to $N=2$ with (\ref{FabEq51}) and
(\ref{FabEq52}), we obtain that uniformly on
$\left\{z:|z-z_k|\geq\epsilon\right\}$ as
$w\to\ou_k$ ($y=w-\ou_k\to 0$)
\begin{equation}\label{FabEq19}
\frac{1}{\psi(w)-z}=\frac{1}{z_k-z}-
\frac{A_ky^{\lambda_k}}{(z_k-z)^2}+\left\{\begin{array}{cc}
\mathcal{O}\left(y^{2\lambda_k}\right),\ &
0<\lambda_k<1 ,\\
 \mathcal{O}\left(y^{\lambda_k+1}\right),\ &
1<\lambda_k<2\,.
\end{array}\right.
\end{equation}
The asymptotic expansion of $\psi'$ about $\ou_k$
is obtained from that of $\psi$ by termwise
differentiation, so that from (\ref{FabEq51}),
(\ref{FabEq6}) and (\ref{FabEq19}) we see that if
$0<\lambda_k<1$, then uniformly on
$\left\{z:|z-z_k|\geq\epsilon\right\}$ as
$w\to\ou_k$,
\begin{equation}\label{FabEq8}
\frac{\psi'(w)}{\psi(w)-z}=\frac{A_k\lambda_ky^{\lambda_k-1}}{z_k-z}+
\frac{\lambda_k\left[2c^{k}_{020}(z_k-z)-A_k^2\right]y^{2\lambda_k-1}}{(z_k-z)^2}+
\mathcal{O}\left(y^{2\lambda_k+\upsilon-1}\right),
\end{equation}
and more specifically,  for $\lambda_k=1/2$,
\begin{eqnarray}\label{FabEq7}
\frac{\psi'(w)}{\psi(w)-z}&=&\frac{A_k\lambda_ky^{\lambda_k-1}}{z_k-z}+
\frac{\lambda_k\left[2c^{k}_{020}(z_k-z)-A_k^2\right]y^{2\lambda_k-1}}{(z_k-z)^2}\nonumber\\
&&+\frac{c^{k}_{111}(1+\lambda_k)y^{\lambda_k}(\log
y)}{z_k-z}+\mathcal{O}\left(y^{\lambda_k}\right).
\end{eqnarray}
Similarly, if $1<\lambda_k<2$, then uniformly on
$\left\{z:|z-z_k|\geq\epsilon\right\}$ as
$w\to\ou_k$
\begin{equation}\label{FabEq11}
\frac{\psi'(w)}{\psi(w)-z}=\frac{A_k\lambda_ky^{\lambda_k-1}}{z_k-z}+
\frac{c^{k}_{110}(1+\lambda_k)y^{\lambda_k}}{z_k-z}+\mathcal{O}\left(y^{2\lambda_k-1}\right).
\end{equation}

For the analytic functions $(w-\ou_k)^\beta$,
$\log(w-\ou_k)$ in $D_\delta(\ou_k)\cap\Delta_1$
corresponding to the branch of the argument
\[\Theta_k-\pi<\arg(w-\ou_k)<\Theta_k+\pi,\quad
w\in\mathbb{C}\setminus\{t\ou_k:t\leq 1\},
\]
let us denote by $(w-\ou_k)^\beta_\pm$ and
$\log_\pm(w-\ou_k)$ their analytic continuations
from $D_\delta(\ou_k)\cap\Delta_1$ onto
$D^\pm_{\delta}(\ou_k)$, respectively. If
$n,\,m\geq 0$ are integers and $\beta>-1$, then
\begin{eqnarray}\label{FabEq9}
&&\hspace{-.5cm}\int_{\sigma_k}^{\ou_k}t^n(t-\ou_k)_\pm^\beta(\log_\pm(t-\ou_k))^mdt\nonumber\\
&=&\int_{0}^{\ou_k}t^n(t-\ou_k)_\pm^\beta(\log_\pm(t-\ou_k))^mdt-\int_{0}^{\sigma_k}t^n(t-\ou_k)_\pm^\beta(\log_\pm(t-\ou_k))^mdt\nonumber\\
&=&e^{\mp i\beta\pi}
e^{i(n+1+\beta)\Theta_k}\int_{0}^{1}x^n(1-x)^\beta
\left[\log(1-x)+i(\Theta_k\mp\pi)\right]^mdx +\mathcal{O}(\sigma^n)\nonumber\\
&=&e^{\mp i\beta\pi}
e^{i(n+1+\beta)\Theta_k}\sum_{l=0}^m{m\choose
l}\alpha_{\beta,m-l}(n)(i(\Theta_k\mp\pi))^l+\mathcal{O}(\sigma^n),
\end{eqnarray}
and
\begin{eqnarray}\label{FabEq18}
\int_{\sigma_k}^{\ou_k}\mathcal{O}\left(t^n(t-\ou_k)_\pm^\beta(\log_\pm(t-\ou_k))^m\right)dt
&=&\mathcal{O}\left(\alpha_{\beta,m}(n)\right)\,.
\end{eqnarray}

Then, we get by combining (\ref{FabEq8}),
(\ref{FabEq9}), (\ref{FabEq18}), (\ref{FabEq16})
and the well-known identity
\[\Gamma(1-z)\Gamma(z)=-z\Gamma(-z)\Gamma(z)=\pi/\sin(\pi z),\]
that if $0<\lambda_k<1$, $\lambda_k\not=1/2$, then
\begin{eqnarray}\label{FabEq23}
&&\hspace{-.5cm}\frac{1}{2\pi
i}\int_{\sigma_k}^{\ou_k}\left(\frac{t^n\psi_+'(t)}{\psi_+(t)-z}-\frac{t^n\psi_-'(t)}{\psi_-(t)-z}\right)dt\nonumber\\
&=& \frac{A_k
\lambda_k\sin(\lambda_k\pi)e^{i(n+\lambda_k)\Theta_k}\alpha_{\lambda_k-1,0}(n)}{\pi(z_k-z)}+\mathcal{O}\left(
\alpha_{2\lambda_k+\upsilon-1,0}(n)\right)
\nonumber\\
&&
+\frac{\left[2c^{k}_{020}(z_k-z)-A_k^2\right]\lambda_k\sin(2\pi\lambda_k)
e^{i(n+2\lambda_k)\Theta_k}\alpha_{2\lambda_k-1,0}(n)}{\pi(z_k-z)^2}\,, \nonumber\\
&=&
\frac{\alpha_{\lambda_k-1,0}(n)}{\Gamma(-\lambda_k)\Gamma(\lambda_k)}\left(\frac{
A_k e^{i(n+\lambda_k)\Theta_k} }{z-z_k}
\right.+o\left(n^{-\lambda_k}\right)\nonumber\\
&&\left.\hspace{3cm}+\frac{\Gamma(-\lambda_k)\left[2c^{k}_{020}(z-z_k)+(A_k)^2\right]e^{i(n+2\lambda_k)\Theta_k}}
{2\Gamma(-2\lambda_k)(z-z_k)^2n^{\lambda_k}}\right).
\end{eqnarray}
Similarly, we get from (\ref{FabEq7}), that if $\lambda_k=1/2$, then
\begin{eqnarray}\label{FabEq47}
&&\hspace{-.5cm}\frac{1}{2\pi
i}\int_{\sigma_k}^{\ou_k}\left(\frac{t^n\psi_+'(t)}{\psi_+(t)-z}-\frac{t^n\psi_-'(t)}{\psi_-(t)-z}\right)dt
\nonumber\\
&=&\frac{\alpha_{\lambda_k-1,0}(n)}{\Gamma(-\lambda_k)\Gamma(\lambda_k)}\left(
\frac{A_k e^{i(n+\lambda_k)\Theta_k}
}{z-z_k}+\frac{c^{k}_{111}(1+\lambda_k)e^{i(n+\lambda_k+1)\Theta_k}(\log
n)} {(z-z_k)n}+o\left(n^{-1}\log n\right)\right),
\end{eqnarray}
and from (\ref{FabEq11}) that if $1<\lambda_k<2$, then
\begin{eqnarray}\label{FabEq58}
&&\hspace{-.5cm}\frac{1}{2\pi
i}\int_{\sigma_k}^{\ou_k}\left(\frac{t^n\psi_+'(t)}{\psi_+(t)-z}-\frac{t^n\psi_-'(t)}{\psi_-(t)-z}\right)dt
\nonumber\\
&=&
\frac{\alpha_{\lambda_k-1,0}(n)}{\Gamma(-\lambda_k)\Gamma(\lambda_k)}\left(
\frac{A_ke^{i(n+\lambda_k)\Theta_k} }{z-z_k}
-\frac{c^{k}_{110}(1+\lambda_k)
e^{i(n+\lambda_k+1)\Theta_k}}
{(z-z_k)n}+o\left(n^{-1}\right)\right).
\end{eqnarray}
Thus, Lemma \ref{FabLem3} for a relevant $\ou_k$ with
$\lambda_k\not\in\{1,2\}$ follows from (\ref{FabEq23}),
(\ref{FabEq47}) and (\ref{FabEq58}).

Next, let us consider the case
$\lambda_k\in\{1,2\}$. From (\ref{FabEq10}) and
(\ref{FabEq13}) we see that \emph{if $m_k\geq
2$}, then uniformly on
$\{z:|z-z_k|\geq\epsilon\}$ as $w\to \ou_k$,

\begin{eqnarray}\label{FabEq15}
\frac{\psi'(w)}{\psi(w)-z}&=&\frac{Q_k'(w)}{Q_k(w)-z}+\frac{\psi'(w)-Q_k'(w)}{Q_k(w)-z}+\mathcal{O}\left(y^{\Lambda_k+\lambda_k-1}(\log
y)^{m_k}\right)\nonumber\\
&=&\frac{Q_k'(w)}{Q_k(w)-z}+\left[\psi'(w)-Q_k'(w)\right]\left(\frac{1}{z_k-z}+\frac{z_k-Q_k(w)}{(z_k-z)(Q_k(w)-z)}\right)\nonumber\\
&&+\mathcal{O}\left(y^{\Lambda_k+\lambda_k-1}(\log
y)^{m_k}\right)\nonumber\\
&=&\frac{Q_k'(w)}{Q_k(w)-z}+\frac{A_{k}\Lambda_k
y^{\Lambda_k-1}(\log
y)^{m_k}}{z_k-z}+\frac{\left(A_k m_k+B_k \Lambda_k\right)y^{\Lambda_k-1}(\log y)^{m_k-1}}{z_k-z}\nonumber\\
&&+\frac{\left[B_k(m_k-1)+C_k\Lambda_k\right]
y^{\Lambda_k-1}(\log
y)^{m_k-2}}{z_k-z}\nonumber\\&&+\left\{\begin{array}{ll}
                \mathcal{O}\left(y^{\Lambda_k-1}(\log y)^{m_k-3}\right) &\quad \mathrm{if}\ m_k\geq 3,\\
               \mathcal{O}\left( y^{\Lambda_k-1/2} \right)& \quad \mathrm{if}\ m_k=
               2.
              \end{array}\right.
\end{eqnarray}
Similarly, one gets from (\ref{FabEq17}) and (\ref{FabEq12}) that
\emph{if $m_k=1$ and $\lambda_k=2$}, then
\begin{eqnarray}\label{FabEq56}
\frac{\psi'(w)}{\psi(w)-z}&=&\frac{Q_k'(w)}{Q_k(w)-z}+\frac{A_{k}\Lambda_k
y^{\Lambda_k-1}(\log
y)^{m_k}}{z_k-z}+\frac{\left(A_k m_k+B_k \Lambda_k\right)y^{\Lambda_k-1}(\log y)^{m_k-1}}{z_k-z}\nonumber\\
&&+\frac{\widetilde{C}_k\left(\Lambda_k+1\right)
y^{\Lambda_k}(\log y)^{\lfloor(r_k+1)/\lambda_k
\rfloor}}{z_k-z}\nonumber\\&&+
\frac{\left[\widetilde{C}_k\lfloor(r_k+1)/\lambda_k
\rfloor+\widetilde{D}_k
\left(\Lambda_k+1\right)\right]y^{\Lambda_k}(\log
y)^{\lfloor(r_k+1)/\lambda_k\rfloor-1}}{z_k-z}\nonumber\\
&&+\left\{\begin{array}{ll}
                          \mathcal{O}\left( y^{\Lambda_k+1/2}\right) & \quad\mathrm{if}\
r_k=2,\\
\mathcal{O}\left( y^{\Lambda_k}(\log y)^{\lfloor(r_k+1)/\lambda_k
\rfloor-2}\right) & \quad\mathrm{if}\ r_k\geq 3,
                        \end{array}\right.
\end{eqnarray}
while \emph{if $m_k=1$ and $\lambda_k=1$}, then
uniformly on $\{z:|z-z_k|\geq\epsilon\}$ as $w\to
\ou_k$,
\begin{eqnarray*}
\frac{1}{\psi(w)-z}&=&\frac{1}{Q_k(w)-z}+\frac{Q_k(w)-\psi(w)}{\left[Q_k(w)-z\right]^2}+
\frac{\left[Q_k(w)-\psi(w)\right]^2}{\left[Q_k(w)-z\right]^2\left[\psi(w)-z\right]}\\
&=&\frac{1}{Q_k(w)-z}+\frac{Q_k(w)-\psi(w)}{\left(z_k-z\right)^2}+\mathcal{O}\left(y^{\Lambda_k}\right),
\end{eqnarray*}
and so
\begin{eqnarray}\label{FabEq20}
\frac{\psi'(w)}{\psi(w)-z}&=&\frac{Q_k'(w)}{Q_k(w)-z}+\left[\psi'(w)-Q_k'(w)\right]\left(\frac{1}{z_k-z}+
\frac{z_k-Q_k(w)}{(z_k-z)^2}\right)\nonumber\\
&&-\frac{c^{k}_{010}A_ky^{\Lambda_k}(\log
y)}{\left(z_k-z\right)^2}+\mathcal{O}
\left(y^{\Lambda_k}\right)\nonumber\\
&=&\frac{Q_k'(w)}{Q_k(w)-z}+\frac{A_{k}\Lambda_ky^{\Lambda_k-1}(\log
y)}{z_k-z}+\frac{\left(A_k+B_k\Lambda_k
\right)y^{\Lambda_k-1}}{z_k-z}\nonumber\\
&&+\frac{\widetilde{C}_k(\Lambda_k+1)y^{\Lambda_k}(\log
y)^{r_k+1}}{z_k-z}\nonumber\\
&&+\frac{\left[\widetilde{C}_k(r_k+1)+ \widetilde{D}_k(\Lambda_k+1)
\right]y^{\Lambda_k}(\log
y)^{r_k}}{z_k-z}\nonumber\\
&&-\frac{c^{k}_{010}A_k\Lambda_ky^{\Lambda_k}(\log
y)}{z_k-z}-\frac{c^{k}_{010}A_ky^{\Lambda_k}(\log
y)}{\left(z_k-z\right)^2}+\mathcal{O}\left(y^{\Lambda_k}(\log
y)^{r_k-1}\right).
\end{eqnarray}

Thus, we get from (\ref{FabEq15}),
(\ref{FabEq56}), (\ref{FabEq20}), (\ref{FabEq9})
and (\ref{FabEq18}) that if $\lambda_k\in\{1,2\}$
and $m_k\geq 2$, then uniformly on
$\{z:|z-z_k|\geq\epsilon\}$ as $w\to \ou_k$
\begin{eqnarray}\label{FabEq26}
&&\hspace{-.5cm}\frac{1}{2\pi
i}\int_{\sigma_k}^{\ou_k}\left(\frac{t^n\psi_+'(t)}{\psi_+(t)-z}
-\frac{t^n\psi_-'(t)}{\psi_-(t)-z}\right)dt\nonumber\\
&=&\frac{m_k\Lambda_k\alpha_{\Lambda_k-1,m_k-1}(n)}{(-1)^{\Lambda_k-1}}\left\{\frac{A_ke^{i(n+\Lambda_k)\Theta_k}}{z-z_k}\right.
+o\left(1/\log n\right)\nonumber\\
&&\hspace{3.6cm}-\left.\frac{(m_k-1)\left(A_km_k+B_k\Lambda_k+i\Theta_km_k\Lambda_kA_k\right)e^{i(n+\Lambda_k)\Theta_k}
}{ m_k\Lambda_k(z-z_k)(\log n)}\right\},
\end{eqnarray}
while if $\lambda_k\in\{1,2\}$ and $m_k=1$, then
uniformly on $\{z:|z-z_k|\geq\epsilon\}$ as $w\to
\ou_k$
\begin{eqnarray}\label{FabEq25}
&&\hspace{-.5cm}\frac{1}{2\pi
i}\int_{\sigma_k}^{\ou_k}\left(\frac{t^n\psi_+'(t)}{\psi_+(t)-z}
-\frac{t^n\psi_-'(t)}{\psi_-(t)-z}\right)dt\nonumber\\
&=&\frac{m_k\Lambda_k\alpha_{\Lambda_k-1,m_k-1}(n)}{(-1)^{\Lambda_k-1}}\left\{\frac{A_ke^{i(n+\Lambda_k)\Theta_k}}{z-z_k}\right.
+ o\left(n^{-1}(\log n)^{\lfloor
(r_k+1)/\lambda_k\rfloor-1}\right)\nonumber\\
&&\hspace{3.5cm}-\left.\frac{\lfloor
(r_k+1)/\lambda_k\rfloor\Gamma(\Lambda_k)\widetilde{C}_k(\Lambda_k+1)e^{i(n+\Lambda_k+1)\Theta_k}}
{m_k\Lambda_k(z-z_k)n(-\log n)^{1-\lfloor
(r_k+1)/\lambda_k\rfloor}}\right\}.
\end{eqnarray}
 This completes the proof of part (a) of Lemma
\ref{FabLem3}. The proof of part (b) easily follows from
(\ref{FabEq53}) by proceeding similarly as in the proof of
part (a).
\end{demo}

\begin{demo}{\emph{\textbf{Proof of Theorem \ref{FabThm1}}}} Let
$E\subset G$ be a compact set, and let
$0<\epsilon<\mathrm{dist}(E,\{z_1,\ldots,z_s\})$.
For this $\epsilon$, choose $\delta$ such that
(\ref{FabEq70}) holds, and fix $\sigma$ with
$0<1-\sigma<\delta$ and satisfying Lemma
\ref{FabLem1}(b).

By (\ref{FabEq27}), if $R>1$, then
\begin{equation}\label{FabEq21}
F_n(z)=\frac{1}{2\pi
i}\ointcc_{\mathbb{T}_R}\frac{t^n\psi'(t)dt}{\psi(t)-z}\quad
\forall\,z\in E.
\end{equation}
Since $\psi'(w)/(\psi(w)-z)$ is analytic on
$\ext(\Gamma_\sigma)$ with continuous boundary
values on
$\Gamma_\sigma\setminus\{\ou_1,\ldots,\ou_s\}$,
and since by (\ref{FabEq1}) and (\ref{FabEq2}),
for $k=1,2,\ldots,s$,
$\psi'(w)=\mathcal{O}\left((w-\ou_k)^{\lambda_k-1}\right)$
as $w\to\ou_k$, then
\begin{equation}\label{FabEq34}
\ointcc_{\mathbb{T}_R}\frac{t^n\psi'(t)dt}{\psi(t)-z}=
\ointcc_{\mathbb{T}_\sigma}\frac{t^n\psi'(t)dt}{\psi(t)-z}+\sum_{k=1}^s
\int_{\sigma_k}^{\ou_k}\left(\frac{t^n\psi_+'(t)}{\psi_+(t)-z}-\frac{t^n\psi_-'(t)}{\psi_-(t)-z}\right)dt,
\end{equation}
so that by Lemma \ref{FabLem3},
\begin{eqnarray}\label{FabEq33}
F_n(z)&=&\mathcal{O}\left(\sigma^n\right)+\sum_{k=1}^s\frac{1}{2\pi
i}
\int_{\sigma_k}^{\ou_k}\left(\frac{t^n\psi_+'(t)}{\psi_+(t)-z}-\frac{t^n\psi_-'(t)}{\psi_-(t)-z}\right)dt
\nonumber\\
&=&\alpha_{\Lambda_1-1,M_1}(n)
\left(\mathcal{C}_1\sum_{k=1}^u\frac{
A_ke^{i(n+\Lambda_k)\Theta_k}}{z-z_k}
+R_n(z)\right),
\end{eqnarray}
where (for every $\tau>0$),
\begin{equation}\label{FabEq82}
R_n(z)=\sum_{k=1}^ur_{\sigma_k,n}(z)+\sum_{k=u+1}^v\mathcal{O}\left(n^{\Lambda_1-\Lambda_k}(\log
n)^{M_k-M_1}\right)+\mathcal{O}\left(n^{-\tau}\right)
\end{equation}
uniformly in $z\in E$ as $n\to\infty$. Theorem
\ref{FabThm1} is proven.
\end{demo}

\begin{demo}{\emph{\textbf{Proof of Theorem \ref{FabThm2}}}}
First, observe that it suffices to prove Theorem
\ref{FabThm2}(a) assuming that $E$ does not
contain $\infty$, because by the very definition
of the Faber polynomials, $F_n(z)-[\phi(z)]^n$ is
analytic at $\infty$, and an application of the
maximum principle for analytic functions will
extend the validity of the theorem to closed sets
of $\Omega\cup L_1$ containing $\infty$.

Then, let $E$ be a closed set ($\infty\not\in E$)
such that either $E\subset \Omega\cup L_1$ or
$E\subset L_2$, and let
\[
0<\epsilon<\mathrm{dist}\left(E,\left\{z_k:z_k\not\in
E \right\}\right).
\]
For this $\epsilon$, choose $\delta$ such that
Lemma  \ref{FabLem2} holds, and fix $\sigma$ with
$0<1-\sigma<\delta$ such that Lemma
\ref{FabLem1}(c) holds.

Recall that  $\eta(z)=1$ if $z\in \Omega \cup
L_1$, $\eta(z)=2$ if $z\in
L_2\setminus\{z_k:\lambda_k=2\}$, and that
$\phi_1(z):=\phi(z)$ for all $z\in \Omega$.

For every $z\in E\setminus\{z_1,\ldots,z_s\}$,
$\phi_l(z)\in \mathrm{ext}(\Gamma_\sigma)$, $l\in
\{1,2,\ldots,\eta(z)\}$, and by Lemma
\ref{FabLem1}(c), the function
$w^n\psi'(w)/\left[\psi(w)-z\right]$ is analytic
in the variable $w$ on
$\mathrm{ext}(\Gamma_\sigma)\setminus\left\{\phi_l(z):1\leq
l\leq \eta(z)\right\}$, with residue
$[\phi_l(z)]^n$ at each (simple pole) $\phi_l(z)$
and continuous boundary values on
$\Gamma_\sigma\setminus\{\ou_1,\ldots,\ou_s\}$.
Moreover, by (\ref{FabEq1}) and (\ref{FabEq2}),
for $k=1,2,\ldots,s$,
$\psi'(w)=\mathcal{O}\left((w-\ou_k)^{\lambda_k-1}\right)$
as $w\to\ou_k$, so that (\ref{FabEq27}) and the
residue theorem yield that \emph{for every $z\in
E\setminus\{z_1,\ldots,z_s\}$},
\begin{equation}\label{FabEq40}
F_n(z)=\sum_{l=1}^{\eta(z)}[\phi_l(z)]^n+\frac{1}{2\pi
i}\ointcc_{\mathbb{T}_\sigma}\frac{t^n\psi'(t)dt}{\psi(t)-z}+\sum_{k=1}^s
\frac{1}{2\pi
i}\int_{\sigma_k}^{\ou_k}\left(\frac{t^n\psi_+'(t)}{\psi_+(t)-z}-\frac{t^n\psi_-'(t)}{\psi_-(t)-z}\right)dt.
\end{equation}

In fact, we claim that \emph{for every $z\in E$},
\begin{eqnarray}\label{FabEq42}
F_n(z)&=&\Phi_n(z)+\frac{1}{2\pi
i}\ointcc_{\mathbb{T}_\sigma}\frac{t^n\psi'(t)dt}{\psi(t)-z}+\sum_{k\,:\,z_k\not\in
E}\int_{\sigma_k}^{\ou_k}\left(\frac{t^n\psi_+'(t)}{\psi_+(t)-z}-\frac{t^n\psi_-'(t)}{\psi_-(t)-z}\right)dt\nonumber\\
&& +\sum_{k\,:\,z_k\in
E}\left(\ointcc_{[\sigma_k,\ou_k]}\frac{t^n
[\psi'(t)-Q'_k(t)]dt}{Q_k(t)-z}
-\ointcc_{[\sigma_k,\ou_k]}\frac{t^n[\psi(t)-Q_k(t)]\psi'(t)dt}{\left(Q_k(t)-z\right)^2}\right.\nonumber\\
&&\hspace{2cm}\left.+\ointcc_{[\sigma_k,\ou_k]}\frac{t^n[\psi(t)-Q_k(t)]^2\psi'(t)dt}
{\left(Q_k(t)-z\right)^2(\psi(t)-z)}\right).
\end{eqnarray}

Indeed, for $z\in E\setminus\{z_1,\ldots,z_s\}$,
the claim is  a direct consequence of
(\ref{FabEq40}) and  the identity
\begin{eqnarray}\label{FabEq29}
\frac{\psi'(w)}{\psi(w)-z}&=&\frac{Q_k'(w)}{Q_k(w)-z}+\frac{\psi'(w)-Q'_k(w)}{Q_k(w)-z}
-\frac{[\psi(w)-Q_k(w)]\psi'(w)}{\left(Q_k(w)-z\right)^2}\nonumber\\
&&+\frac{[\psi(w)-Q_k(w)]^2\psi'(w)}{\left(Q_k(w)-z\right)^2(\psi(w)-z)},
\end{eqnarray}
taking into account that, by (\ref{FabEq94}) and
(\ref{FabEq31}), $Q_k(w)-z\not=0$ for all $w\in
[\sigma_k,\ou_k)$, $z\in \cj{\Omega}$

Suppose now $j$ is such that $z_j\in E$, and let
us agree in that, in case $\eta(z_j)=2$ and one
of the two values $\phi_1(z_j),\phi_{2}(z_j)$ is
not contained in $\{\ou_1,\ldots,\ou_s\}$, that
value is precisely $\phi_1(z_j)$.

Then, for every $k$ with $z_k=\psi(\ou_k)=z_j$
(there are at most two of them), let $T_k\subset
D_{\delta}(\ou_k)$ be the circle centered at
$\ou_k$ of radius $1-\sigma$. We can assume
$\sigma$ was chosen so close to 1 that in case
$z_j\in L_2$ and $\phi_1(z_j)\not\in
\{\ou_1,\ldots,\ou_s\}$, $\phi_1(z_j)$ lies in
the exterior of $T_k$.

Then, we obtain once again from (\ref{FabEq27})
and Lemma \ref{FabLem1}(c) that
\begin{eqnarray}\label{FabEq57}
F_n(z_j)&=&\frac{1}{2\pi
i}\ointcc_{\mathbb{T}_\sigma}\frac{t^n\psi'(t)dt}{\psi(t)-z_j}+\sum_{k\,:\,z_k\not=z_j}
\frac{1}{2\pi
i}\int_{\sigma_k}^{\ou_k}\left(\frac{t^n\psi_+'(t)}{\psi_+(t)-z_j}-\frac{t^n\psi_-'(t)}{\psi_-(t)-z_j}\right)dt
\nonumber\\
&&+\sum_{k\,:\,z_k=z_j}\frac{1}{2\pi
i}\ointcc_{T_k}\frac{t^n\psi'(t)dt}{\psi(t)-z_k}+
\left\{\begin{array}{ll} [\phi_1(z_j)]^n, &\
\mathrm{if}\ \eta(z_j)=2,\ \phi_1(z_j)\not\in\{\ou_1,\ldots,\ou_s\},\\
0,&\ \mathrm{otherwise}.
\end{array}
\right.\nonumber\\
\end{eqnarray}

By Lemma \ref{FabLem2}(a),
$w^nQ'_k(w)/\left(Q_k(w)-z_k\right)$ is analytic
on $D^*_{\delta}(\ou_k)$ with a simple pole at
$\ou_k$, so that if we take identity
(\ref{FabEq29}) for $z=z_k$, multiply it by $w^n$
and integrating it over $T_k$, we obtain that for
$k$ with $z_k=z_j$,
\begin{eqnarray}\label{FabEq89}
\ointcc_{T_k}\frac{t^n\psi'(t)dt}{\psi(t)-z_j}&=&2\pi
i\lambda_k(\ou_k)^n+\ointcc_{[\sigma_k,\ou_k]}\frac{t^n
[\psi'(t)-Q'_k(t)]dt}{Q_k(t)-z_j}
-\ointcc_{[\sigma_k,\ou_k]}\frac{t^n[\psi(t)-Q_k(t)]\psi'(t)dt}{\left(Q_k(t)-z_j\right)^2}\nonumber\\
&&+\ointcc_{[\sigma_k,\ou_k]}\frac{t^n[\psi(t)-Q_k(t)]^2\psi'(t)dt}
{\left(Q_k(t)-z_j\right)^2(\psi(t)-z_j)}\,.
\end{eqnarray}
Then, (\ref{FabEq42}) for $z=z_j$ follows from
relations (\ref{FabEq57}) and (\ref{FabEq89}).

Now that the claim is proven, we proceed to estimate the
integrals that occur in (\ref{FabEq42}) under the symbol
$\sum_{k\,:\,z_k\in E}$ . For this, we first observe that
if $\{\mathcal{F}(\cdot,z):z\in E\}$ is a uniformly bounded
family of measurable functions on $[\sigma_k,\ou_k)$, if
$n,\,m\geq 0$ are integers and $\beta>-1$, then (compare to
(\ref{FabEq9}))
\begin{eqnarray}\label{FabEq36}
&&\hspace{-.5cm}\int_{\sigma_k}^{\ou_k}\mathcal{F}(t,z)t^n(t-\ou_k)_\pm^\beta(\log_\pm(t-\ou_k))^mdt\nonumber\\
&=&e^{\mp i\beta\pi}
e^{i(n+1+\beta)\Theta_k}\sum_{l=0}^m{m\choose
l}(i(\Theta_k\mp\pi))^l\alpha_{\beta,m-l}(n)G_{\beta,l,n}(z)+\mathcal{O}(\sigma^n),
\end{eqnarray}
where the functions
\[
G_{\beta,l,n}(z):=\frac{\int_{0}^1\widetilde{\mathcal{F}}(x,z)x^n(1-x)^\beta(\log
(1-x))^{m-l}dx}{\alpha_{\beta,m-l}(n)},\quad
\widetilde{\mathcal{F}}(x,z):=\left\{\begin{array}{ll}
                                        \mathcal{F}\left(xe^{i\Theta_k},z\right), & \ x\in [\sigma,1),\\
                                       1, &  \ x\in [0,\sigma),
                                     \end{array}
\right.
\]
are uniformly bounded on $E$, are independent of
the sign $\pm$, and $G_{\beta,l,n}(z)=1$ whenever
$\mathcal{F}(\cdot,z)\equiv 1$.

Now, assume $\ou_k$ is relevant. Recall that with
$y=w-\ou_k$ (see (\ref{FabEq14}), (\ref{FabEq17}),
(\ref{FabEq13}) and (\ref{FabEq12})),
\begin{equation}\label{FabEq44}
Q_k(w)=z_k+c^{k}_{010}y^{\lambda_k}+\mathcal{O}\left(y^{\lambda_k+1}\right),\quad
Q_k'(w)=c^{k}_{010}\lambda_k
y^{\lambda_k-1}+\mathcal{O}\left(y^{\lambda_k}\right),
\end{equation}
\begin{eqnarray}\label{FabEq43}
\psi(w)-Q_k(w)&=&A_{k}y^{\Lambda_k}(\log
y)^{m_k}+B_{k}y^{\Lambda_k}(\log
y)^{m_k-1}\nonumber\\
&&+\left\{\begin{array}{cc}
\mathcal{O}\left(y^{\Lambda_k}(\log y)^{m_k-2}\right)\ &\ \mathrm{if}\ m_k\geq 2,\\
\mathcal{O}\left(y^{\Lambda_k+1}(\log y)^{\lfloor
(r_k+1)/\lambda_k\rfloor}\right)\ &\ \mathrm{if}\
m_k=1,
\end{array}\right.
\end{eqnarray}
\begin{eqnarray}\label{FabEq41}
\psi'(w)-Q'_k(w)&=&A_{k} \Lambda_k
y^{\Lambda_k-1}(\log
y)^{m_k}+D_{k}y^{\Lambda_k-1}(\log
y)^{m_k-1}\nonumber\\
&&+\left\{\begin{array}{cc}
\mathcal{O}\left(y^{\Lambda_k-1}(\log y)^{m_k-2}\right)\ &\ \mathrm{if}\ m_k\geq 2,\\
\mathcal{O}\left(y^{\Lambda_k}(\log y)^{\lfloor
(r_k+1)/\lambda_k\rfloor}\right)\ &\ \mathrm{if}\
m_k=1,
\end{array}\right.
\end{eqnarray}
where $c^{k}_{010}\not=0$, $A_k\not=0$, $B_k$ and
$D_k$ are certain constants.

If we set
$\mathcal{F}(t,z):=(Q_k(t)-z_k)/(Q_k(t)-z)$, then
by (\ref{FabEq31}) in Lemma \ref{FabLem2},
(\ref{FabEq44}), (\ref{FabEq41}) and the equality
\[
\left(Q_k(t)-z_k\right)^{-1}=\left(c^{k}_{010}\right)^{-1}(t-\ou_k)^{-\lambda_k}+\mathcal{O}\left((t-\ou_k)^{1-\lambda_k}\right),\qquad
t\to\ou_k\,,
\]
we have that
\begin{eqnarray*}
\int_{\sigma_k}^{\ou_k}\frac{t^n[\psi_\pm'(t)-Q'_k(t)]dt}{Q_k(t)-z}&=&
\int_{\sigma_k}^{\ou_k}\frac{\mathcal{F}(t,z)t^n[\psi_\pm'(t)-Q'_k(t)]dt}{Q_k(t)-z_k}\nonumber\\
&=&\frac{A_{k}
\Lambda_k}{c^{k}_{010}}\int_{\sigma_k}^{\ou_k}\mathcal{F}(t,z)t^n(t-\ou_k)_\pm^{r_k-1}(\log_\pm
(t-\ou_k))^{m_k}dt\nonumber\\
&&+\frac{D_{k}}{c^{k}_{010}}\int_{\sigma_k}^{\ou_k}\mathcal{F}(t,z)t^n(t-\ou_k)_\pm^{r_k-1}(\log_\pm
(t-\ou_k))^{m_k-1}dt\nonumber\\
&&+\left\{\begin{array}{cc}
\int_{\sigma_k}^{\ou_k}\mathcal{O}\left(t^n(t-\ou_k)_\pm^{r_k-1}(\log_\pm (t-\ou_k))^{m_k-2}\right),\ &\ m_k\geq 2, \\
\int_{\sigma_k}^{\ou_k}\mathcal{O}\left(t^n(t-\ou_k)_\pm^{r_k}(\log_\pm
(t-\ou_k))^{\lfloor
(r_k+1)/\lambda_k\rfloor}\right),\ &\ m_k=1.
\end{array}\right.
\end{eqnarray*}
Combining this with (\ref{FabEq36}) we see that
uniformly in $z\in E$ as $n\to\infty$,
\begin{equation}\label{FabEq39}
\ointcc_{[\sigma_k,\ou_k]}\frac{t^n
[\psi'(t)-Q'_k(t)]dt}{Q_k(t)-z}=2\pi
i\alpha_{r_k-1,m_k-1}(n)\left( \frac{A_{k}
\Lambda_k m_k
e^{i(n+r_k)\Theta_k}G_{r_k-1,1,n}(z)}{(-1)^{r_k-1}c^{k}_{010}}+o(1)\right),
\end{equation}
where $G_{r_k-1,1,n}(z)=1$ if $z=z_k$.

Similarly, we get from (\ref{FabEq43}), (\ref{FabEq41}) and the
equality
\[
\left(Q_k(t)-z_k\right)^{-2}=\left(c^{k}_{010}\right)^{-2}(t-\ou_k)^{-2\lambda_k}+\mathcal{O}(t-\ou_k),\qquad
t\to\ou_k,
\]
that
\begin{eqnarray*}
\int_{\sigma_k}^{\ou_k}\frac{t^n[\psi(t)_\pm-Q_k(t)]\psi_\pm'(t)dt}{\left(Q_k(t)-z\right)^2}&=&
\frac{A_{k}
\lambda_k}{c^{k}_{010}}\int_{\sigma_k}^{\ou_k}\left[\mathcal{F}(t,z)\right]^2t^n(t-\ou_k)_\pm^{r_k-1}(\log_\pm
(t-\ou_k))^{m_k}dt\\
&&+\frac{\lambda_kB_{k}}{c^{k}_{010}}\int_{\sigma_k}^{\ou_k}\left[\mathcal{F}(t,z)\right]^2t^n(t-\ou_k)_\pm^{r_k-1}(\log_\pm
(t-\ou_k))^{m_k-1}dt\\
&&+\left\{\begin{array}{cc}
\int_{\sigma_k}^{\ou_k}\mathcal{O}\left(t^n(t-\ou_k)_\pm^{r_k-1}(\log_\pm (t-\ou_k))^{m_k-2}\right),\ &\ m_k\geq 2, \\
\int_{\sigma_k}^{\ou_k}\mathcal{O}\left(t^n(t-\ou_k)_\pm^{r_k}(\log_\pm
(t-\ou_k))^{2m_k}\right),\ &\ m_k=1.
\end{array}\right.
\end{eqnarray*}
Hence, uniformly in $z\in E$ as $n\to\infty$,
\begin{equation}\label{FabEq37}
\ointcc_{[\sigma_k,\ou_k]}\frac{t^n[\psi(t)_\pm-Q_k(t)]\psi_\pm'(t)dt}{\left(Q_k(t)-z\right)^2}=2\pi
i\alpha_{r_k-1,m_k-1}(n)\left( \frac{A_{k}
\lambda_k m_k
e^{i(n+r_k)\Theta_k}G_{r_k-1,1,n}(z)}{(-1)^{r_k-1}c^{k}_{010}}+o(1)\right)
\end{equation}
where $G_{r_k-1,1,n}(z)=1$ if $z=z_k$.

As for the last integral in (\ref{FabEq42}), it follows directly
from (\ref{FabEq43}), (\ref{FabEq41}), (\ref{FabEq31}) and
(\ref{FabEq36}) that
\begin{eqnarray}\label{FabEq45}
\int_{\sigma_k}^{\ou_k}\frac{t^n[\psi_\pm(t)-Q_k(t)]^2\psi_\pm'(t)dt}{\left(Q_k(t)-z\right)^2(\psi(t)_\pm-z)}
&=&\int_{\sigma_k}^{\ou_k}\left(\frac{\psi_\pm(t)-z_k}{\psi_\pm(t)-z}\right)
\frac{\left[\mathcal{F}(t,z)\right]^2[\psi_\pm(t)-Q_k(t)]^2\psi_\pm'(t)dt}{\left(Q_k(t)-z_k\right)^2(\psi_\pm(t)-z_k)}
\nonumber\\&=&\int_{\sigma_k}^{\ou_k}
\mathcal{O}\left(t^n(t-\ou_k)_\pm^{2r_k-1}(\log_\pm(t-\ou_k))^{2m_k}\right)dt
\nonumber\\&=&\mathcal{O}\left(\alpha_{2r_k-1,2m_k}(n)\right)\,.
\end{eqnarray}

If $\ou_k$  is not relevant, the degree of $Q_k(z)$  may be
assumed to be as large as desired (see paragraph following
(\ref{FabEq53})), and a similar (easier) analysis shows
that the integrals in the left-hand sides of
(\ref{FabEq39}), (\ref{FabEq37}), (\ref{FabEq45}) are
$\mathcal{O}\left(n^{-\tau}\right)$ uniformly in $z\in E$
as $n\to\infty$, where $\tau$ can be taken arbitrarily
large.

With this last observation in mind, we then obtain by
combining (\ref{FabEq42}), Lemma \ref{FabLem3},
(\ref{FabEq39}), (\ref{FabEq37}) and (\ref{FabEq45}) that
\begin{eqnarray}\label{FabEq24}
F_n(z)&=&\Phi_n(z)+\mathcal{O}(\sigma^n)+\sum_{1\leq
k\leq v\,:\,z_k\not\in
E}\alpha_{\Lambda_k-1,M_k}(n)
\left(\sum_{k=1}^u\frac{\mathcal{C}_kA_k
e^{i(n+\Lambda_k)\Theta_k}}{z-z_k} +o(1)\right)\nonumber\\
&&+ \sum_{1\leq k\leq v\,:\,z_k\in
E}\alpha_{r_k-1,M_k}(n)\left( \frac{A_{k} r_k m_k
e^{i(n+r_k)\Theta_k}G_{r_k-1,1,n}(z)}{(-1)^{r_k-1}c^{k}_{010}}+o(1)\right)
+\mathcal{O}\left(n^{-\tau}\right)
\end{eqnarray}
uniformly on $z\in E$ as $n\to\infty$, where
$\tau$ can be taken arbitrarily large, the
functions $G_{r_k-1,1,n}(z)$ are uniformly
bounded on $E$ and $G_{r_k-1,1,n}(z_k)=1$.

Taking into account that
$(r_k,M_k)<(\Lambda_k,M_k)$, it is now clear that
Theorem \ref{FabThm2}(a) follows from
(\ref{FabEq24}).

It only remains to prove part (b) of Theorem \ref{FabThm2}.
Let $\epsilon>0$ be such that
\[
D_{\epsilon}(z_k)\cap
D_{\epsilon}(z_j)=\emptyset\quad\mathrm{ whenever}\quad
z_k\not=z_j.
\]
For this $\epsilon$, choose $\delta>0$ such Lemma
\ref{FabLem2} holds, and choose $\sigma$ with
$0<1-\sigma<\delta$ such that Lemma \ref{FabLem1}(c) holds.

By increasing $\sigma$ toward $1$ if necessary, we can
assume that $\Gamma_\sigma$ is such that for every
$k=1,2,\ldots,s$, the elements of
$\left\{\phi_1(z_k),\ldots,\phi_{\eta(z_k)}(z_k)\right\}\setminus\{\ou_k\}$
lie outside the circle $T_k\subset D_{\delta}(\ou_k)$
centered at $\ou_k$ with radius $1-\sigma$.

Now, think of $z_j$ as being fixed, so that by Lemma
\ref{FabLem1}(c), the function
\[w^n\psi'(w)/\left(\psi(w)-z_j\right)\] is analytic in the
variable $w$ on
$\mathrm{ext}(\Gamma_\sigma)\setminus\left\{\phi_1(z_j),\ldots,\phi_{\eta(z_j)}(z_j)\right\}$
with a simple pole at each $\phi_l(z_j)\not\in
\{\ou_1,\ldots,\ou_s\}$ and residue $[\phi_l(z_j)]^n$.
Hence, we obtain from (\ref{FabEq27}) that
\begin{eqnarray}\label{FabEq46}
F_n(z_j)&=&
\sum_{\phi_l(z_j)\not\in\{\ou_1,\ldots,\ou_s\}}[\phi_l(z_j)]^n+
\sum_{k\,:\,z_k\not=z_j}\frac{1}{2\pi
i}\int_{\sigma_k}^{\ou_k}\left(\frac{t^n\psi_+'(t)}{\psi_+(t)-z_j}-\frac{t^n\psi_-'(t)}{\psi_-(t)-z_j}\right)dt
\nonumber\\&& +\sum_{k\,:\,z_k=z_j}\frac{1}{2\pi
i}\ointcc_{T_k}\frac{t^n\psi'(t)dt}{\psi(t)-z_j}+\mathcal{O}(\sigma^n).
\end{eqnarray}

If $k$ is such that $z_k=z_j$ and $\lambda_k\in \{1,2\}$,
then as we have previously seen,
\begin{eqnarray}\label{FabEq38}
\frac{1}{2\pi
i}\ointcc_{T_k}\frac{t^n\psi'(t)dt}{\psi(t)-z_j}&=&\lambda_k(\ou_k)^n+\left\{
\begin{array}{ll}
{\displaystyle\alpha_{r_k-1,M_k}(n)\left(
\frac{A_{k} r_k m_k
e^{i(n+r_k)\Theta_k}}{(-1)^{r_k-1}c^{k}_{010}}+o(1)\right)},
&\ \mathrm{if}\ \ou_k\ \mathrm{is\ relevant},\\
  \mathcal{O}\left(n^{-\tau}\right), &\ \mathrm{otherwise},
\end{array}\right.\nonumber\\
\end{eqnarray}
as $n\to\infty$.

If $k$ is such that $z_k=z_j$ and  $\lambda_k\not \in
\{1,2\}$, then we have in virtue of (\ref{FabEq51}) and
(\ref{FabEq52}) that for all $\upsilon>0$ small enough,
\begin{equation*}
\psi(w)=z_j+A_k
y^{\lambda_k}+c^{k}_{020}y^{2\lambda_k}+
\mathcal{O}\left(y^{2\lambda_k+\upsilon}\right),\quad
\mathrm{if}\ \,0<\lambda_k<1,
\end{equation*}
\begin{equation*}
\psi(w)=z_j+A_k
y^{\lambda_k}+c^{k}_{110}y^{1+\lambda_k}+c^{k}_{020}y^{2\lambda_k}+\mathcal{O}\left(y^{2\lambda_k+\upsilon}\right),\quad
\mathrm{if}\ \,1<\lambda_k<2.
\end{equation*}
Hence, if
$0<\upsilon<\min\{\lambda_k,2-\lambda_k\}\ (<1)$,
then
\begin{equation*}
\frac{\psi'(w)}{\psi(w)-z_j}=\frac{\lambda_k}{y}+\frac{c^{k}_{020}\lambda_k
y^{\lambda_k-1}}{A_k}+\mathcal{O}\left(y^{\lambda_k+\upsilon-1}\right),\quad
\mathrm{if}\ \,0<\lambda_k<1,
\end{equation*}
\begin{equation*}
\frac{\psi'(w)}{\psi(w)-z_j}=\frac{\lambda_k}{y}+\frac{c^{k}_{110}}{A_k}+\frac{c^{k}_{020}\lambda_k
y^{\lambda_k-1}}{A_k}+\mathcal{O}\left(y^{\lambda_k+\upsilon-1}\right),\quad
\mathrm{if}\ \,1<\lambda_k<2,
\end{equation*}
so that
\begin{eqnarray}\label{FabEq72}
\frac{1}{2\pi
i}\ointcc_{T_k}\frac{t^n\psi'(t)dt}{\psi(t)-z_j}&=&\lambda_k\left(\ou_k\right)^n+\frac{c^{k}_{020}\lambda_k
}{2\pi i A_k}\ointcc_{[\sigma_k,\ou_k]}t^n(t-\ou_k)^{\lambda_k-1}dt\nonumber\\
&&+\ointcc_{[\sigma_k,\ou_k]}\mathcal{O}\left(t^n(t-\ou_k)^{\lambda_k+\upsilon-1}\right)dt\nonumber\\
&=&\lambda_k\left(\ou_k\right)^n-\alpha_{\lambda_k-1,0}(n)\left(\frac{c^{k}_{020}e^{i(n+\lambda_k)\Theta_k}
}{\Gamma(-\lambda_k)A_k}+o(1)\right).
\end{eqnarray}
Then, coupling  Lemma \ref{FabLem3} with (\ref{FabEq46}),
(\ref{FabEq38}) and (\ref{FabEq72}) yields for every
$\tau>0$
\begin{eqnarray}\label{FabEq55}
F_n(z_j)&=& \Phi_n(z_j) +\sum_{1\leq k\leq
v\,:\,z_k\not=z_j}\alpha_{\Lambda_k-1,M_k}(n)
\left(\frac{\mathcal{C}_kA_k
e^{i(n+\Lambda_k)\Theta_k}}{z_j-z_k} +o(1)\right)
\nonumber\\&& +\sum_{1\leq k\leq v\,:\,z_k=z_j,\
\lambda_k\in\{1,2\}}\alpha_{r_k-1,M_k}(n)\left( \frac{A_{k}
r_k m_k
e^{i(n+r_k)\Theta_k}}{(-1)^{r_k-1}c^{k}_{010}}+o(1)\right)\nonumber\\&&
-\sum_{k\,:\,z_k=z_j,\
\lambda_k\not\in\{1,2\}}\alpha_{\Lambda_k-1,M_k}(n)\left(\frac{c^{k}_{020}e^{i(n+\lambda_k)\Theta_k}
}{\Gamma(-\lambda_k)A_k}+o(1)\right)+\mathcal{O}\left(n^{-\tau}\right)+\mathcal{O}(\sigma^n).
\end{eqnarray}
Theorem \ref{FabThm1}(c) follows immediately by
comparing the terms in (\ref{FabEq55}).
\end{demo}

\begin{demo}{\emph{\textbf{Proof of Theorem
\ref{FabThm8}}}} Given a closed set $E\subset
\left(\Omega\cup L_1\right)\setminus
\{z_1,\ldots, z_s\}$, ($\infty\not\in
E$)\footnote{See remark at the beginning of the
proof of Theorem \ref{FabThm2}.} let
\[E_1^{-1}:=\{w\in \mathbb{T}_1:\psi(w)\in
E\}.\] Then, to prove the theorem it suffices to
show that there is an open set $U\supset
E_1^{-1}$ such that $\phi$ has an analytic an
univalent continuation to $\Omega\cup U$ and
formula (\ref{FabEq76}) holds uniformly in $z\in
E\cup U$ as $n\to\infty$.

The set $\left\{w\in \mathbb{T}_1:\psi(w)\in
L_1\setminus \{z_1,\ldots,z_s\}\right\}$ is the
union of finitely many pairwise disjoint open
circular arcs $J_1,J_2,\ldots,J_{s'}$. Let
$e_1,e_2,\ldots ,e_{s''}$ be those points of
$\mathbb{T}_1$ that happen to be an endpoint of
some $J_l$. Fix
$0<\epsilon<\mathrm{dist}\left(E,\{z_1,\ldots,z_s\}\right)$,
and for this $\epsilon$ choose a number $\delta$
for which Lemma \ref{FabLem2} holds. We can
assume $\delta$ is so small that the disks
$D_\delta(e_1),\ldots, D_\delta(e_{s''})$ are
pairwise disjoint and for every
$e_k\not\in\{\ou_1,\ldots,\ou_s\}$, the analytic
continuation of $\psi$ to $D_\delta(e_k)$
satisfies
\begin{enumerate}
\item[i)]
$\cj{\psi\left(D_\delta(e_k)\right)}\subset
\mathbb{C}\setminus E$.
\end{enumerate}

Let $J_l^*:=J_l\setminus \cup_{k=1}^{s''}D_\delta(e_k)$
($1\leq l\leq s'$). Choose $\sigma$  such that
$1-\delta<\sigma<1/\sigma<1+\delta$ and Lemma \ref{FabLem1}
holds. Let
\[
O_{\delta,\sigma}:=\left\{w:\sigma<|w|<1/\sigma,\quad
w\not\in \cup_{l=1}^{s''}
\cj{D_\delta(e_l)},\quad w\not\in \cup_{k=1}^{s}
\cj{D_\delta(\ou_k)}\right\}
\]
and let $O^l_{\delta,\sigma}$ ($1\leq l\leq s'$)
be  the unique open component of
$O_{\delta,\sigma}$ that contains points of the
arc $J_l^*$. Since $\psi(J_l^*)\subset L_1$,
$\psi\left(\mathbb{T}_1\setminus J_l^*\right)\cap
\psi\left(J_l^*\right)=\emptyset $ for all $1\leq
l\leq s'$, and therefore, if $\sigma $ was taken
sufficiently close to $1$, then for $1\leq l\leq
s'$,
\begin{enumerate}
\item[ii)]
$\psi\left(O^l_{\delta,\sigma}\right)
\cap\psi\left(O_{\delta,\sigma}\setminus
O^l_{\delta,\sigma}\right)=\emptyset$ and
$\psi\left(O^l_{\delta,\sigma}\cap\mathbb{D}_1\right)\subset
G$ .
\end{enumerate}

It follows from ii) that $\psi$ is univalent on
$\Delta_1\cup_{k=1}^{s''}O^l_{\delta,\sigma}$,
and consequently, $\phi$ has an analytic and
univalent continuation to $\Omega\cup\psi
\left(\cup_{k=1}^{s''}O^l_{\delta,\sigma}\right)$.
Let $V\subset
\cup_{l=1}^{s''}O^l_{\delta,\sigma}$ be an open
set containing $E_1^{-1}$, such that $U:=\psi(V)$
intercepts neither
$\cup_{k=1}^{s}\cj{D_\epsilon(\ou_k)}$ nor
$\cup_{l=1}^{s''}\cj{\psi(D_\delta(e_l))}$.

Then, from Lemma \ref{FabLem1}(a), i) and ii)
above, we see that for every $z\in E\cup U$,
$\phi(z)\in \mathrm{ext}(\Gamma_\sigma)$, and the
function $w^n\psi'(w)/\left[\psi(w)-z\right]$ is
analytic in the variable $w$ on
$\mathrm{ext}(\Gamma_\sigma)\setminus\left\{\phi(z)\right\}$
with a simple pole at $\phi(z)$ and residue
$[\phi(z)]^n$. Moreover, it has continuous
boundary values on
$\Gamma_\sigma\setminus\{\ou_1,\ldots,\ou_s\}$
and, by (\ref{FabEq1}) and (\ref{FabEq2}), for
every $k=1,2,\ldots,s$,
$\psi'(w)=\mathcal{O}\left((w-\ou_k)^{\lambda_k-1}\right)$
as $w\to\ou_k$, so that from (\ref{FabEq27}) and
the residue theorem we obtain
\begin{equation*}
F_n(z)=[\phi(z)]^n+\frac{1}{2\pi
i}\ointcc_{\mathbb{T}_\sigma}\frac{t^n\psi'(t)dt}{\psi(t)-z}+\sum_{k=1}^s
\frac{1}{2\pi
i}\int_{\sigma_k}^{\ou_k}\left(\frac{t^n\psi_+'(t)}{\psi_+(t)-z}-\frac{t^n\psi_-'(t)}{\psi_-(t)-z}\right)dt
\end{equation*}
uniformly in $z\in E\cup U$ as $n\to\infty$. The
theorem follows from this and Lemma
\ref{FabLem3}.
\end{demo}

\begin{demo}{\emph{\textbf{Proof of Theorem
\ref{FabThm9}}}} Let $E\subset L_2\setminus
\{z_1,\ldots,z_s\}$ be an analytic arc and let
$0<\epsilon<\mathrm{dist}\left(E,\{z_1,\ldots,z_s\}\right)$.
For this $\epsilon$, choose $\delta>0$ such that Lemma
\ref{FabLem2} holds. Note that $E$ is a subarc of a larger
analytic arc $F\subset L_2$ with all points of $E$ being
interior points of $F$. Notice also that the set
$\{w\in\mathbb{T}_1:\psi(w)\in E\}$ consists of two
disjoint closed circular arcs contained in
$\mathbb{T}_1\setminus\{\ou_1,\ldots,\ou_s\}$, that we
denote by $E^{-1}_1$, $E^{-1}_2$.

By the Schwarz reflection principle, we can find
disjoint open sets $V_1\supset E_1^{-1}$,
$V_2\supset E_2^{-1}$, such that $\psi$ has an
analytic continuation from $\Delta_1$ to
$\Delta_1\cup V_1\cup V_2$, which satisfies
\begin{enumerate}
\item[(i)] $\psi$ is univalent on each $V_i$, $i=1,2$,
 \begin{equation}\label{FabEq96}
\psi(V_1)=\psi(V_2)\subset \mathbb{C}\setminus \cup_{k=1}^s
D_{\epsilon}(z_k),
 \end{equation}
and
\[
\psi(V_1\cap \mathbb{D}_1)=\psi(V_2\cap \Delta_1)\subset
\Omega, \quad \psi(V_1\cap \Delta_1)=\psi(V_2\cap
\mathbb{D}_1)\subset \Omega.
\]
\end{enumerate}
Hence,
 \begin{equation}\label{FabEq98}
\psi\left(\Delta_1\setminus (V_1\cup V_2)\right)\cap
\psi(V_1)=\emptyset.
 \end{equation}

Also, if $Y:=\mathbb{T}_1\setminus \left(V_1\cup
V_2\cup_{k=1}^s D_\delta(\ou_k)\right)$, then $\psi(Y)$ is
a closed set disjoint from $E$, and we can find a
neighborhood $W_Y$ of $Y$ and open sets $V_1'$, $V_2'$

\begin{enumerate}
\item[(ii)] $Y\subset W_Y$, $E_1^{-1}\subset V_1'\subset V_1$, $E_2^{-1}\subset V_2'\subset
V_2$, $\psi(V'_1)=\psi(V'_2)$, $(V'_1\cup
V'_2)\cap Y=\emptyset$, and the analytic
continuation of $\psi$ from $\Delta_1$ to
$\Delta_1\cup W_Y$ satisfies
\begin{equation}\label{FabEq97}
\psi(W_Y)\cap\psi(V'_1\cup V'_2)=\emptyset.
\end{equation}
\end{enumerate}

Then, the sets $W_Y$, $V_1$, $V_2$, $D_{\delta}(\ou_1),
D_{\delta}(\ou_2),\ldots, D_{\delta}(\ou_s)$, form an open
covering of $\mathbb{T}_1$, and we can find $\sigma$ with
$0<1-\sigma<\delta$, such that Lemma \ref{FabLem1}(a) holds
and
\[
\{w:\sigma\leq |w|<1\}\subset W_Y\cup V_1\cup
V_2\cup\left(\cup_{k=1}^sD_{\delta}(\ou_k)\right).
\]
Let $V''_1$ and $V''_2$ be two open sets of the form
\[
V''_1=\{\rho_1<|w|<1/\rho_1\,,\  \alpha_1<
\arg(w)<\beta_1\},
\]
\[
V''_2=\{\rho_2<|w|<1/\rho_2\,,\  \alpha_2<
\arg(w)<\beta_2\},
\]
such that
\[
E_1^{-1}\subset V''_1\subset V'_1\cap
\ext(\Gamma_\sigma),\quad E_2^{-1}\subset V''_2\subset
V'_2\cap \ext(\Gamma_\sigma),
\]
and let $U:=\psi(V''_1)\cap\psi(V''_2)$, so that $E\subset
U\subset \psi(V_1)$.

If $z\in U$ and $\psi(w)=z$ for some
$w\in\Delta_\sigma$, then by (\ref{FabEq70}) and
(\ref{FabEq96}) we must have $w\in
W_Y\cup\left(\Delta_1\setminus (V_1\cup
V_2)\right)\cup V_1\cup V_2$. But, since
$U\subset\psi(V'_1)\subset \psi(V_1)$,
(\ref{FabEq98}) and (\ref{FabEq97}) imply that
indeed $w\in V_1\cup V_2$. Since $\psi$ is
one-to-one on each $V_i$, $i=1,2$, and
$V''_1\subset V_1\cap \ext(\Gamma_\sigma)$,
$V''_2\subset V_1\cap \ext(\Gamma_\sigma)$, we
conclude that $w$ must belong to $V''_1\cup
V''_2$. Hence, if $z\in U$ and $w_{z,1}\in
V''_1$, $w_{z,2}\in V''_2$ are such that
$\psi(w_{z,1})=\psi(w_{z,2})=z$, then the
function $(\psi(w)-z)^{-1}$ is analytic on
$\ext(\Gamma_\sigma)\setminus\{w_{z,1},w_{z,2}\}$,
with simple poles at $w_{z,1},w_{z,2}$ and
continuous boundary values on $\Gamma_\sigma$.
Then again, from (\ref{FabEq27}) and the residue
theorem we obtain
\begin{equation}\label{FabEq92}
F_n(z)=(w_{z,1})^n+(w_{z,1})^n+\frac{1}{2\pi
i}\ointcc_{\mathbb{T}_\sigma}\frac{t^n\psi'(t)dt}{\psi(t)-z}+\sum_{k=1}^s
\frac{1}{2\pi
i}\int_{\sigma_k}^{\ou_k}\left(\frac{t^n\psi_+'(t)}{\psi_+(t)-z}-\frac{t^n\psi_-'(t)}{\psi_-(t)-z}\right)dt
\end{equation}
uniformly in $z\in U$ as $n\to\infty$.

Now, let $\phi_+$, $\phi_-$ be, respectively, the inverse
functions of $\psi|_{V''_1}$, $\psi|_{V''_2}$. Then both
$\phi_+$, $\phi_-$ are defined on $U$ and for $z\in U$,
$\phi_+(z)=w_{z,1}$, $\phi_-(z)=w_{z,2}$. But, by the
Schwarz reflection principle, if
\[
U^+:=\psi(V''_1\cap \Delta_1)\cap\psi(V''_2\cap
\mathbb{D}_1)\subset \Omega,\quad  U^-:=\psi(V''_2\cap
\Delta_1)\cap\psi(V''_1\cap \mathbb{D}_1)\subset \Omega,
\]
then $U^+$ is the reflection of $U^-$ about $E$ and
$\phi_{\pm}$ is the analytic continuation of
$\phi|_{U^\pm}$ across $E$ to all of $U$.  The theorem
follows from (\ref{FabEq92}) and Lemma \ref{FabLem3}.
\end{demo}

\begin{demo}{\emph{\textbf{Proof of Theorem \ref{FabThm3}}}}

In this particular case, $\psi(w)=(w^s+1)^{1/s}$
and
\[
\ou_k=e^{i(2k-1)\pi /s}, \quad z_k=0,
\quad\lambda_k=1/s\quad \forall\,
k\in\{1,\ldots,s\}.
\]

Let the compact set $E\subset G=\{z:|z^s-1|<1\}$
be given, fix $\epsilon$ with
$0<\epsilon<\mathrm{dist}(E,\{z_1,\ldots,z_s\})$,
and for this $\epsilon$ choose $\delta$
satisfying Lemma \ref{FabLem2}. Then, as shown in
the proof of Theorem \ref{FabThm1}, for every
$1-\delta<\sigma<1$ sufficiently close to 1, we
have
\begin{equation}\label{FabEq71}
F_n(z)=\mathcal{O}(\sigma^n)+\sum_{k=1}^s
\frac{1}{2\pi
i}\int_{\sigma_k}^{\ou_k}\left(\frac{t^n\psi_+'(t)}{\psi_+(t)-z}-\frac{t^n\psi_-'(t)}{\psi_-(t)-z}\right)dt
\end{equation}
uniformly in $z\in E$ as $n\to\infty$.

Notice that $\psi(w)=e^{-2\pi i k
/s}\psi\left(we^{2\pi i k /s}\right)$ for
$k=1,\ldots,s$. Then, from identity
(\ref{FabEq5}) if $N=2s-l$ we get that for
$n=sm+l$ ($1\leq l\leq s-1$), then
\[
\int_{\sigma_k}^{\ou_k}\frac{t^n\psi'_\pm(t)dt}{\psi_\pm(t)-z}=-\sum_{j=1}^{2s-l}z^{-j}\int_{\sigma_k}^{\ou_k}
t^n\psi'_\pm(t)\left[\psi_\pm(t)\right]^{j-1}dt
+\int_{\sigma_k}^{\ou_k}\frac{t^n\psi'_\pm(t)\left[\psi_\pm(t)\right]^{2s-l}dt}{z^{2s-l}\left[\psi_\pm(t)-z\right]}.
\]
Hence,
\begin{eqnarray}\label{FabEq73}
\sum_{k=1}^{s}\int_{\sigma_k}^{\ou_k}\frac{t^n\psi'_\pm(t)dt}{\psi_\pm(t)-z}&=&-\sum_{j=1}^{2s-l}\left(\sum_{k=1}^{s}e^{2\pi
\,
i(k-1)(n+j)/s}\right)z^{-j}\int_{\sigma_1}^{\ou_1}
t^n\psi'_\pm(t)\left[\psi_\pm(t)\right]^{j-1}dt
\nonumber\\
&&+\int_{\sigma_1}^{\ou_1}\mathcal{O}\left(t^n(t-\ou_1)^{(s-l+1)/s}\right)dt\nonumber\\
&=&-z^{l-s}\int_{\sigma_1}^{\ou_1}
t^n\psi'_\pm(t)\left[\psi_\pm(t)\right]^{s-l-1}dt+\mathcal{O}\left(\alpha_{(s+1-l)/s,0}(n)\right)\nonumber\\
&&-z^{l-2s}\int_{\sigma_1}^{\ou_1}
t^n\psi'_\pm(t)\left[\psi_\pm(t)\right]^{2s-l-1}dt.
\end{eqnarray}
Now,
\[
(w^s+1)^{1/s}=s^{1/s}e^{i\pi(s-1)/s^2}(w-\ou_1)^{1/s}\left[1+\frac{(s-1)e^{-i\pi/s}(w-\ou_1)}{2s}+\mathcal{O}\left((w-\ou_1)^2\right)\right],
\]
and therefore,
\begin{eqnarray}\label{FabEq74}
\psi'(w)[\psi(w)]^{s-l-1}&=&s^{-l/s}e^{i\pi(s-1)(s-l)/s^2}\left[(w-\ou_1)^{-l/s}
+\frac{(s-1)(2s-l)(w-\ou_1)^{(s-l)/s}}{2se^{i\pi/s}}\right.\nonumber\\
&&\left.
\hspace{3cm}+\mathcal{O}\left((w-\ou_1)^{(2s-l)/s}\right)
\right]
\end{eqnarray}
and
\begin{equation}\label{FabEq75}
\psi'(w)[\psi(w)]^{2s-l-1}=s^{(s-l)/s}e^{i\pi(s-1)(2s-l)/s^2}(w-\ou_1)^{(s-l)/s}+\mathcal{O}\left((w-\ou_1)^{(2s-l)/s}\right).
\end{equation}

Theorem \ref{FabThm3} is then obtained after somewhat
lengthy computations by combining (\ref{FabEq71}),
(\ref{FabEq73}), (\ref{FabEq74}), (\ref{FabEq75}),
(\ref{FabEq9}) and (\ref{FabEq18}).
\end{demo}

We conclude this section with the proof of
(\ref{FabEq16}).

\begin{demo}{\emph{\textbf{Proof of (\ref{FabEq16})}}} We shall prove the equivalent
statement that for any two integers $n,\,m\geq 0$
and real $\beta>-1$,
\begin{equation}\label{FabEq48}
\alpha_{\beta,m}(n)=\frac{\Gamma(\beta+1)n!(-\log
(n+\beta+1))^m\left[1+\varepsilon_{\beta,m}(n)/\log(n+\beta+1)\right]}{\Gamma(n+\beta+2)},
\end{equation}
where $\varepsilon_{\beta,m}(n)=\mathcal{O}(1)$
as $n\to\infty$. The proof is by induction on
$m$. For $m=0$, (\ref{FabEq48}) is trivially
satisfied with $\varepsilon_{\beta,m}(n)\equiv
0$. Assume that (\ref{FabEq48}) is also true for
some $m\geq 0$. Integrating by parts we obtain
\begin{equation*}
\alpha_{\beta,m+1}(n)=\frac{n\alpha_{\beta+1,m+1}(n-1)}{\beta+1}-\frac{(m+1)\alpha_{\beta,m}(n)}{\beta+1}.
\end{equation*}
Since,
$\alpha_{\beta+1,m+1}(n-1)=\alpha_{\beta,m+1}(n-1)-\alpha_{\beta,m+1}(n)$,
this gives
\begin{equation}\label{FabEq4}
\alpha_{\beta,m+1}(n)=\frac{n\alpha_{\beta,m+1}(n-1)}{n+\beta+1}-\frac{(m+1)\alpha_{\beta,m}(n)}{n+\beta+1}.
\end{equation}
By iterating (\ref{FabEq4}) and from the
induction hypothesis we then obtain
\begin{eqnarray}\label{FabEq49}
\alpha_{\beta,m+1}(n)&=&\frac{\Gamma(\beta+2)n!\,\alpha_{\beta,m+1}(0)}{\Gamma(n+\beta+2)}
-(m+1)\sum_{j=1}^n\frac{\Gamma(i+\beta+1)n!\,\alpha_{\beta,m}(i)}{\Gamma(n+\beta+2)i!}\nonumber\\
&=&\frac{\Gamma(\beta+1)n!(m+1)}{\Gamma(n+\beta+2)(-1)^{m+1}}\left(\frac{m!}{(\beta+1)^{m+1}}
+\sum_{j=1}^n\frac{\left(\log(j+\beta+1)\right)^m}{j+\beta+1}\right.\nonumber\\
&&\hspace{4.5cm}\left.+\sum_{j=1}^n\frac{\varepsilon_{\beta,m}(j)\left(\log(j+\beta+1)\right)^{m-1}}{j+\beta+1}\right).
\end{eqnarray}
Then, the validity of (\ref{FabEq48}) follows
immediately from (\ref{FabEq49}) and the fact
that for all $m\geq 0$,
\begin{eqnarray*}
\sum_{j=1}^n\frac{(\log(j+\beta+1))^{m}}{j+\beta+1}&=&\int_{1}^{n}\frac{\left(\log
(x+\beta+1)\right)^{m}dx}{x+\beta+1}+\mathcal{O}(1)\\
&=&
\frac{(\log(n+\beta+1))^{m+1}}{m+1}+\mathcal{O}(1)
\end{eqnarray*}
as $n\to\infty$.
\end{demo}

\section{Proofs of the zero
results}\label{FabSec3}

\begin{demo}{\emph{\textbf{Proof of Corollary
\ref{FabPro1}}}} Suppose there is a compact set
$E\subset G$ and a subsequence
$\{n_j\}\subset\mathbb{N}$ such that
$F^*_{n_j}(z)$ has more that $J-1$ zeros on $E$
counting multiplicities, where $J$ is the number
of corners of $L$ ($J\leq s$). By assumption A.3
(and extracting a subsequence from $\{n_j\}$ if
needed), we can assume that $\{F^*_{n_j}\}$
converges locally uniformly on $G$ to a nonzero
rational function $R(z)$ with denominator having
degree no larger than $J-1$. By Hurwitz's
Theorem, there is an open set $U\supset E$ such
that for all $j$ large enough, $F^*_{n_j}$ and
$R(z)$ have the same number of zeros on $U$,
contradicting our assumption.

We now show that $\nu_{n}\wc \mu$, for which we
use standard arguments. By Helly's selection
theorem \cite[Thm. 1.3]{Saff}, from every
subsequence of $\{\nu_{n}\}_{n\geq 1}$ it is
possible to extract another subsequence
converging in the weak*-topology to a measure
$\mu$. Thus, to finish the proof, it suffices to
show that every such limit measure $\mu$ is the
equilibrium measure $\mu_L$ of $L$.

Then, suppose $\nu_{n_j}\wc \mu$ as $j\to\infty$,
so that by Corollary \ref{FabCor2} and what we
just proved above, $\mu$ must be supported on
$L$. Let us denote by $U^\alpha(z)$ the
logarithmic potential of the measure $\alpha$,
that is,
\[
U^{\alpha}(z):=\int\log
\frac{1}{|z-t|}d\alpha(t).
\]
Then, if $\kappa_n$ denotes the leading
coefficient of $F_n$, then
$\kappa_n=[\phi'(\infty)]^n$, and we obtain from
(\ref{FabEq68}) and the fact that $\nu_{n_j}\wc
\mu$, that for all $z\in\Omega$
\[
U^{\mu}(z)=\lim_{j\to\infty}U^{\nu_{n_j}}(z)
=\lim_{j\to\infty}\frac{1}{n_j}\log\frac{\kappa_{n_j}}{|F_{n_j}(z)|}=
\log|\phi'(\infty)/\phi(z)|\,.
\]
On the other hand, it is not difficult to see
from the definition of $\mu_L$ in (\ref{FabEq77})
that for all $z\in \Omega$,
$U^{\mu_L}(z)=\log|\phi'(\infty)/\phi(z)|$.
Hence, $\mu$ and $\mu_L$ are two measures
supported on $L$ whose logarithmic potential
coincide on $\Omega$, which in view of Carleson's
theorem \cite[Thm. 4.13]{Saff} forces
$\mu=\mu_L$.
\end{demo}

\begin{demo}{\emph{\textbf{Proof of Theorem \ref{FabThm4}}} }
Suppose that for some subsequence
$\{n_\nu\}_{\nu\geq 1}\subset \mathbb{N}$,
\[
H_{n_\nu}(z)=\sum_{k=1}^u \frac{\hat{A}_k e^{2\pi
in_\nu\theta_k}}{z-z_k}\to f(z)\quad
\mathrm{as}\quad \nu\to\infty.
\]
By extracting a subsequence if necessary, we may
assume that for some fixed $\ell\in
\{1,\ldots,\mathbf{q}\}$,
$n_\nu=\mathrm{\mathbf{q}}m_\nu+\ell$ with
$m_\nu\in \mathbb{N}$, and by the compacity of
$\mathbb{T}_1$, that for some real numbers
$\alpha_2,\ldots,\alpha_{u^*}$
\begin{equation*}
\lim_{\nu\to\infty}e^{2\pi i r_{kj}n_{\nu}
\theta_j}=e^{2\pi i r_{kj}\alpha_j},\quad 1\leq
k\leq u,\ 2\leq j\leq u^*,
\end{equation*}
so that by (\ref{FabEq99}), $f$ must have the
form (\ref{FabEq100}).

Conversely, we now show that given an integer
$\ell\in \{1,\ldots,\mathbf{q}\}$ and arbitrary
real numbers $\alpha_2,\ldots,\alpha_{u^*}$, it
is possible to choose a subsequence
$\{n_\nu\}_{\nu\geq 1}$ such that
\begin{equation}\label{FabEq101}
\lim_{\nu\to\infty}e^{2\pi i n_\nu
\theta_k}=e^{2\pi i\,\left(\frac{n_\nu
p_k}{q_k}+\sum_{j=2}^{u^*}r_{kj}n_\nu\theta_j\right)}=e^{2\pi
i\,\left(\frac{\ell
p_k}{q_k}+\sum_{j=2}^{u^*}r_{kj}\alpha_j\right)}\,,\quad
 1\leq k\leq u.
\end{equation}

For this, we first observe that given arbitrary
real numbers $\chi_2,\ldots,\chi_u$, it is always
possible to find a subsequence
$\{m_\nu\}_{\nu\geq 1}\subset\mathbb{N}$ such
that
\begin{equation}\label{FabEq79}
\lim_{\nu\to\infty}e^{2\pi i
r_{kj}\mathrm{\mathbf{q}}m_{\nu}
\theta_j}=e^{2\pi i
r_{kj}\mathrm{\mathbf{q}}\chi_j},\quad 1\leq
k\leq u,\ 2\leq j\leq u^*.
\end{equation}
In effect, consider the set of linear forms in
the variable $x$
\[\left\{r_{kj}\mathrm{\mathbf{q}}\theta_jx: 1\leq
k\leq u,\ 2\leq j\leq u^*\right\},\] and suppose
$\beta_{kj}$, $1\leq k\leq u$, $2\leq j\leq u^*$,
are integers such that
\[
\sum_{k,j}
\beta_{kj}r_{kj}\mathrm{\mathbf{q}}\theta_jx=
x\sum_{j=2}^{u^*}\left(\sum_{k=1}^u\beta_{kj}r_{kj}\mathrm{\mathbf{q}}\right)\theta_j
\]
is a linear form whose coefficient is an integer.
Then, by the linear independence of the numbers
$1,\theta_2,\ldots\theta_{u^*}$, we must have
$\sum_{k=1}^u\beta_{kj}r_{kj}\mathrm{\mathbf{q}}=0$
for every $2\leq j\leq u^*$. Hence, for an
arbitrary collection of real numbers
$\chi_2,\ldots,\chi_{u^*}$, we have $\sum_{k,j}
\beta_{kj}r_{kj}\mathrm{\mathbf{q}}\theta_j\chi_j=0$,
and so by Kronecker's theorem \cite[Chap. III,
Thm. IV.]{Cassels}, it is possible to find a
subsequence  $\{m_\nu\}_{\nu\geq 1}$ satisfying
(\ref{FabEq79}).

Then, choose a subsequence $\{m_\nu\}_{\nu\geq
1}$ satisfying (\ref{FabEq79}) with
$\chi_j=(\alpha_j-\ell\theta_j)/\mathrm{\mathbf{q}}$,
$2\leq j\leq u^*$. Then, (\ref{FabEq101}) is
satisfied by the subsequence
$n_\nu:=\mathrm{\mathbf{q}}m_\nu+\ell$, $\nu\in
\mathbb{N}$.

It only remains to prove that there is a rational
function $f$ of the form (\ref{FabEq100}) that is
not identically zero. Assume without loss of
generality that the set $\{k:z_k=z_1,\ 1\leq
k\leq u\}$ consists of the numbers
$1,2,\ldots,u'$ for some $u'\leq u$. It suffices
to show that it is impossible to have
\begin{equation}\label{FabEq86} \sum_{k=1}^{u'}\hat{A}_k e^{2\pi
i\,\left(\frac{\ell
p_k}{q_k}+\sum_{j=2}^{u^*}r_{kj}\alpha_j\right)}=0,
\quad \ell\in \{1,\ldots,\mathbf{q}\},\
\alpha_2,\ldots,\alpha_{u^*}\in \mathbb{R}.
\end{equation}
Assume, on the contrary, that this is the case.
Since $\hat{A}_1\not=0$ and $r_{1j}=0$ for all
$2\leq j\leq u^*$, we must obviously have
$r_{kj}=0$ for all $1\leq k\leq u'$, $2\leq j\leq
u^*$, and consequently,
\[
\theta_k=\frac{p_k}{q_k}\,,\qquad
k=1,2,\ldots,u'.
\]

Let $q'\leq \mathrm{\mathbf{q}}$ be the least
common multiple of the denominators
$q_1,q_2,\ldots,q_{u'}$, and for
$k=1,2,\ldots,u'$, set $p'_k:=p_k q'/q_k$, so
that $1\leq p'_k\leq q'$, and since
$\theta_1,\theta_2,\ldots,\theta_{u'}$ are
pairwise distinct, so are the numbers
$p'_1,p'_2,\ldots,p'_{u'}$, and therefore $u'\leq
q'$. Then, by (\ref{FabEq86}) we must have
\[
\sum_{k=1}^{u'}\hat{A}_k \left(e^{2\pi i
\ell/q'}\right)^{p'_k}=0 \quad \forall\,
\ell\in\{1,\ldots, q'\}.
\]
But this homogenous system of linear equations
with unknowns $\hat{A}_k$'s has only the trivial
solution, since the Vandermonde matrix
$\left(a_{l,m}\right)_{1\leq\ell,m\leq q'\ }$,
$a_{l,m}=\left(e^{2\pi i\,\ell /q'}\right)^m$, is
nonsingular. This contradicts that all the
$\hat{A}_k $'s are nonzero.
\end{demo}

\begin{demo}{\emph{\textbf{Proof of Corollary \ref{FabCor1}}}}
By Theorem \ref{FabThm4}, there is a subsequence
$\{n_j\}\subset \mathbb{N}$ such that
$\{F^*_{n_j}\}$ converges locally uniformly on
$G$ to a nonzero rational function. Then,
proceeding exactly as in the proof of Corollary
\ref{FabPro1}, we find that $\nu_{n_j}\wc \mu_L$
as $j\to\infty$.
\end{demo}

\begin{demo}{\emph{\textbf{Proof of Corollary \ref{FabCor3}}}
} This is just a straightforward consequence of Theorem
\ref{FabThm4} and Hurwitz's theorem, therefore, we omit it.
\end{demo}

\noindent \textbf{Acknowledgement.} The author
extends his gratitude to Prof. A.
Mart\'{\i}nez-Finkelshtein for valuable
discussions relevant to this work.

\noindent E. Mi\~{n}a-D\'{\i}az\\
Department of Mathematical Sciences\\
Indiana University-Purdue University\\
Fort Wayne, IN 46805\\
USA\\
minae@ipfw.edu
\end{document}